\documentclass[twoside]{article}
\usepackage{amsmath}			%voor wiskundige tekens
\usepackage{amssymb}			%voor wiskundige tekens
\usepackage{amsthm}			%voor theorem environments
\usepackage{a4wide}			%voor uitlijning pagina
\usepackage{enumitem}			%voor opsommingen
\usepackage{dsfont}			%certain fonts
\usepackage{verbatim}			%voor comments
\usepackage[all]{xy}			%voor commutatieve diagrammen
\usepackage{fancyhdr}          	 	%voor opmaak o.a. headers en footnotes
\usepackage[english]{babel}     	%spelling

\numberwithin{equation}{section}

\title{A classification of $SU(d)$-type C$^*$-tensor categories}
\author{Bas Jordans}
\date{\today}

\newtheorem{Thm}{Theorem}[section]
\newtheorem{Prop}[Thm]{Proposition}
\newtheorem{Lemma}[Thm]{Lemma}
\newtheorem{Cor}[Thm]{Corollary}

\theoremstyle{definition}
\newtheorem{Assum}[Thm]{Assumption}
\newtheorem{Def}[Thm]{Definition}
\newtheorem{Exam}[Thm]{Example}
\newtheorem{Rem}[Thm]{Remark}
\newtheorem{Not}[Thm]{Notation}

\newenvironment{pf}{\vspace{0 pt plus 1pt minus 1pt}{\it Proof. }}{\mbox{}\hfill$\boxtimes$\vspace{6pt plus 1pt minus 2 pt}}

\newcommand{\Nat}{\mathbb{N}}
\newcommand{\Int}{\mathbb{Z}}

\newcommand{\Rea}{\mathbb{R}}
\newcommand{\Com}{\mathbb{C}}
\newcommand{\Torus}{\mathbb{T}}

\newcommand{\Ccal}{\mathcal{C}}
\newcommand{\Dcal}{\mathcal{D}}

\newcommand{\Hcal}{\mathcal{H}}

\newcommand{\ra}{\rightarrow}

\newcommand{\ten}{\otimes}
\newcommand{\Cstar}{C$^*$-}
\newcommand{\unit}{\mathds{1}}

\DeclareMathOperator{\End}{End}
\DeclareMathOperator{\im}{Im}
\DeclareMathOperator{\Hom}{Hom}
\DeclareMathOperator{\Ob}{Ob}
\DeclareMathOperator{\ran}{ran}

\DeclareMathOperator{\Rep}{Rep}

\DeclareMathOperator{\Span}{span}
\DeclareMathOperator{\Tr}{Tr}
\DeclareMathOperator{\tr}{tr}

\begin{document}
\thispagestyle{plain}
\begin{center}
{\Large A classification of $SU(d)$-type C$^*$-tensor categories}

\bigskip
{\large Bas Jordans\footnote{E-mail address: bpjordan@math.uio.no\\
Department of Mathematics, University of Oslo, P.O. Box 1053 Blindern, 0316 Oslo, Norway. \\ 
The research leading to these results has received funding from the European Research Council under the European
Union’s Seventh Framework Programme (FP/2007-2013) / ERC Grant Agreement no. 307663 (PI: S. Neshveyev).\\
Date: \today}
}
\end{center}

\bigskip
\begin{abstract}
Kazhdan and Wenzl classified all rigid tensor categories with fusion ring isomorphic to the fusion ring of the group $SU(d)$. In this paper we consider the C$^*$-analogue of this problem. Given a rigid C$^*$-tensor category $\mathcal{C}$ with fusion ring isomorphic to the fusion ring of the group $SU(d)$, we can extract a constant $q$ from $\mathcal{C}$ such that there exists a $*$-representation of the Hecke algebra $H_n(q)$ into $\mathcal{C}$. The categorical trace on $\mathcal{C}$ induces a Markov trace on $H_n(q)$. Using this Markov trace and a representation of $H_n(q)$ in $\textrm{Rep}\,(SU_{\sqrt{q}}(d))$ we show that $\mathcal{C}$ is equivalent to a twist of the category $\textrm{Rep}\,(SU_{\sqrt{q}}(d))$. Furthermore a sufficient condition on a C$^*$-tensor category $\mathcal{C}$ is given for existence of an embedding of a twist of $\textrm{Rep}\,(SU_{\sqrt{q}}(d))$ in $\mathcal{C}$.
\end{abstract}

\section{Introduction}
Tannaka-type reconstruction theorems allow one to reconstruct an algebraic object (for example a group) from its category of representations. There are numerous of these theorems, the classical Tannaka-Krein duality \cite{Tannaka39}, the Doplicher-Roberts theorem \cite{DoplicherRoberts89}, Deligne's theorem \cite{Deligne90}, Woronowicz's duality for compact matrix pseudogroups \cite{Woronowicz88} and many more. 
Despite these theorems it is still very difficult (if not impossible) to give a complete list of all quantum groups which satisfy the fusion rules of a certain group. However, if one instead tries to classify all (\Cstar) tensor categories which have a fusion ring isomorphic to the fusion ring of a certain group $G$, this problem becomes easier to solve. Kazhdan and Wenzl \cite{KazhdanWenzl} gave such a classification in the case of tensor categories with fusion ring isomorphic to the fusion ring $K[\Rep(SU(d))]$. They showed that if $\Ccal$ is a tensor category with fusion ring isomorphic to $K[\Rep(SU(d))]$, then there exists a constant $\mu\in\Com^*$ not a non-trivial root of unity such that $\Ccal$ is (monoidally) equivalent to a ``twist'' of $\Rep(SU_\mu(d))$, the representation category of the quantum group $SU_\mu(d)$. These twists are determined by a $d$-th root of unity. \\
This paper contains two main results. The first one (cf. Theorem \ref{classification}) is the \Cstar analogue of the result by Kazhdan and Wenzl. We will show that all \Cstar tensor categories which satisfy the same fusion rules as $\Rep(SU(d))$, the so-called $SU(d)$-type categories, can be classified by pairs $(\mu,\omega)$ where $\mu\in(0,1]$ and $\omega$ is a $d$-th root of unity. Namely given a $SU(d)$-type category $\Ccal$ we can extract constants $\mu$ and $\omega$ from $\Ccal$ such that $\Ccal$ is equivalent to a ``twist'' by $\omega$ of the category $\Rep(SU_{\mu}(d))$. The other main theorem is inspired on the paper by Pinzari \cite{Pinzari07}. 
In this paper she gives a sufficient condition when it is possible to embed $\Rep(SU_\mu(d))$ in a given braided \Cstar tensor category. We generalize this result to conditions on \Cstar tensor categories which are sufficient to construct an embedding of a ``twist'' of $\Rep(SU_\mu(d))$ in a given category (see Theorem \ref{embedding}). These two main results are independent of each other, but the proofs of both of theorems are related. They are both based on Theorem \ref{characterization} which gives some technical conditions when a given category is equivalent to a ``twist'' of $\Rep(SU_\mu(d))$. In this paper Hecke algebras play a key-role, because the representation category $\Rep(SU_\mu(d))$ has a natural representation of the Hecke algebra. Following Kazhdan and Wenzl we construct a representation of the Hecke algebra in the the endomorphisms of a $SU(d)$-type category. These representations allow us eventually to recover the category from its fusion ring. In these categories we need to make some explicit calculations. However, in general categories this is often very difficult. Therefore we use the categorical trace \cite{LongoRoberts97} to show that the representation of the Hecke algebra is independent of the category $\Ccal$. This result allows us to make computations in $\Rep(SU_\mu(d))$ in which everything is more explicit. \\
This paper is organized as follows. In Section \ref{tensor_cat} we start by recalling the main definitions and properties of \Cstar categories and specialize to $SU(d)$-type categories. We continue with the necessary results on Hecke algebras. Section \ref{computations} will be devoted to making the necessary computations in $\Rep(SU_\mu(d))$ allowing us later to compute the twist invariant of a general $SU(d)$-type category. In section \ref{representation} we construct the representation of the Hecke algebra into $\End(X^{\ten n})$ and we establish that this representation is independent of the category $\Ccal$. In the next section we consider a specific class of \Cstar tensor categories and we prove a technical theorem showing that all these categories of this specific type are equivalent to a ``twist'' of $\Rep(SU_\mu(d))$. Section \ref{2characterizations} contains the two main theorems of this paper.

\section{Preliminaries on C$^*$-tensor categories}\label{tensor_cat}
In this section we will introduce the specific class of \Cstar tensor categories we are interested in, namely the $SU(d)$-type categories. We will also discuss ``twists'' of such categories.

We will not give the full definitions of \Cstar tensor categories and functors of these categories. Precise definitions of \Cstar tensor categories and functors thereof can be found in e.g., \cite[\textsection 2.1]{NeshveyevTuset13}. The essential element one has to keep in mind is that in a \Cstar tensor category one is able to take ``tensor products'' of objects and morphisms and that there exists a conjugation operation. We will define the associativity morphisms and conjugate objects in \Cstar tensor categories, because they will play a key-role later on.  

\begin{Def}\label{Cstar-tensor_category}
A \Cstar tensor category is a \Cstar category equipped with a bilinear bifunctor $\ten \colon \Ccal\times\Ccal\ra\Ccal$, $(U,V)\mapsto U\ten V$, which will be called the {\it tensor product} and it is required that there exist natural unitary isomorphisms
\[
\alpha_{U,V,W}\colon(U\ten V)\ten W \ra U\ten (V\ten W),
\]
the {\it associativity morphisms} such that the {\it pentagonal diagram}
\[
\xymatrix{& ((U\ten V)\ten W)\ten X \ar[ld]_{\alpha\ten\iota} \ar[rd]^{\alpha_{12,3,4}}& \\
(U\ten(V\ten W))\ten X\ar[d]_{\alpha_{1,23,4}} & & (U\ten V)\ten (W\ten X)\ar[d]^{\alpha_{1,2,34}}\\
U\ten ((V\ten W)\ten X)\ar[rr]^{\iota\ten\alpha} & & U\ten (V\ten (W\ten X))
}
\]
commutes. Here the convention of leg-numbering is used (e.g., $\alpha_{12,3,4}:= \alpha_{U\ten V,W,X}$). Furthermore, it is assumed that there exists an object $\unit$ (the {\it unit}) and natural unitary isomorphisms
\[
\lambda_U\colon\unit\ten U\ra U,\qquad \rho_U\colon U\ten\unit\ra U
\]
such that $\lambda_\unit=\rho_\unit$ and the {\it triangle diagram}
\[
\xymatrix{(U\ten\unit)\ten V \ar[rr]^{\alpha}\ar[rd]_{\rho\ten\iota} & & U\ten (\unit\ten V) \ar[ld]^{\iota\ten\lambda}\\
&U\ten V &
}
\]
commutes. A category will be called {\it strict} if
\[
(U\ten V)\ten W = U\ten (V\ten W), \qquad \unit\ten U = U = U\ten\unit
\]
and the associativity morphisms $\alpha$ and the morphisms $\lambda$ and $\rho$ are the identity morphisms. We assume that \Cstar tensor categories are closed under subobjects and direct sums and that the unit object $\unit$ is simple. 
\end{Def}

\begin{Rem}\label{subobjects}
Sometimes we will use the terminology that an object $U$ is a subobject of an object $V$, or simply $U\subset V$. What is meant by this is that there exists a projection $p\in\End(V)$ and a morphism $v\in\Hom(U,V)$ such that $v^*v=id_U$ and $vv^*=p$. Via this $v$ we can restrict morphisms $T\in\End(V)$ to the object $U$, we write $T|_U:=v^*Tv\in\End(U)$.
\end{Rem}

\begin{Def}
Let $\Ccal$ and $\Dcal$ be \Cstar tensor categories. A functor $F:\Ccal\ra\Dcal$ together with an isomorphism $F_0\colon\unit_\Dcal\ra F(\unit_\Ccal)$ and natural isomorphisms $F_2\colon F(U)\ten F(V)\ra F(U\ten V)$ is called a {\it \Cstar tensor functor} if for any $U,\,V,\,W\in\Ob(\Ccal)$, the diagrams 
\[
\xymatrix{(F(U)\ten F(V))\ten F(W)\ar[r]^-{F_2\ten\iota}\ar[d]_{\alpha_\Dcal} & F(U\ten V)\ten F(W)\ar[r]^{F_2} & F((U\ten V)\ten W)\ar[d]^{F(\alpha_\Ccal)}\\
F(U)\ten (F(V)\ten F(W))\ar[r]^-{\iota\ten F_2} & F(U)\ten F(V\ten W)\ar[r]^{F_2} & F(U\ten (V\ten W))
}
\]
\[
\xymatrix{\unit_\Dcal\ten F(U)\ar[r]^{F_0\ten\iota}\ar[d]_{\lambda_\Dcal}& F(\unit_\Ccal)\ten F(U)\ar[d]^{F_2}  & & F(U)\ten \unit_\Dcal\ar[r]^{\iota\ten F_0} \ar[d]_{\rho_\Dcal} & F(U)\ten F(\unit_\Ccal)\ar[d]^{F_2}\\
F(U) & F(\unit_\Ccal\ten U)\ar[l]_{F(\lambda_\Ccal)} & & F(U) & F(U\ten \unit_\Ccal)\ar[l]_{F(\rho_\Ccal)}
}
\]
commute. The \Cstar tensor functor $F$ is called {\it unitary} if $F(T^*)=F(T)^*$ on morphisms and $F_2\colon F(U)\ten F(V)\ra F(U\ten V)$ and $F_0$ are unitary. $F$ is called  {\it fully faithful} if $F\colon\Hom(U,V)\ra \Hom(F(U),F(V))$ is an isomorphism. $F$ is called {\it essentially surjective} if for every object $U\in\Ob(\Dcal)$ there exists an object $V\in\Ob(\Ccal)$ such that $U$ is isomorphic to $F(V)$. $F$ is called a {\it monoidal equivalence} if $F$ is fully faithful and essentially surjective. Two \Cstar tensor categories $\Ccal$ and $\Dcal$ are {\it monoidally equivalent} if there exists a monoidal equivalence $F\colon\Ccal\ra\Dcal$.
\end{Def}

\begin{Rem}
Any \Cstar tensor category can be strictified \cite[\textsection XI.3]{MacLane98}. This means that if $\Ccal$ is a (non-strict) \Cstar tensor category, then there exists a strict \Cstar tensor category $\Dcal$ such that $\Ccal$ and $\Dcal$ are unitarily monoidally equivalent. So unless stated otherwise we deal with strict categories.\\
If $\Ccal$ is a category which satisfies all the requirements of a \Cstar tensor category except from the existence of direct sums and subobjects, then $\Ccal$ can be completed to a new category which is a \Cstar tensor category (see for example, \cite[\textsection 2.5]{NeshveyevTuset13}). For this define $\Ccal'$ with $\Ob(\Ccal'):=\{(U_1,\ldots, U_n)\,:\, n\geq1,\, U_i\in\Ob(\Ccal)\}$ and $\Hom_{\Ccal'}((U_1,\ldots, U_m), (V_1,\ldots, V_n)):= \bigoplus_{i,j}\Hom_\Ccal(U_i,V_j)$. Now $(U_i)_i\oplus (V_j)_j:= (U_1,\ldots, U_m,V_1,\ldots, V_n)$ and $(U_i)_i\ten(V_j)_j$ is given by the tuple consisting of the lexicographical ordering of $U_i\ten V_j$. Let $\Ccal''$ be the category with $\Ob(\Ccal''):=\{(U,p)\,:\, U\in\Ob(\Ccal'),\,p\in\End_{\Ccal'}(U) \textrm{ projection}\}$ and $\Hom_{\Ccal''}((U,p),(V,q)):=q\Hom_{\Ccal'}(U,V)p$. The tensor product of objects is given by $(U,p)\ten(V,q):=(U\ten V, p\ten q)$. The involution, direct sums and tensor products of morphisms on $\Ccal'$ and $\Ccal''$ are defined in the obvious way. Then $\Ccal''$ is a \Cstar tensor category. It is clear that there exists a unitary tensor functor $i:\Ccal\ra\Ccal''$. \\
The completion $\Ccal''$ is universal in the following sense: if $\Dcal$ is a \Cstar tensor category and $F\colon\Ccal\ra\Dcal$ is a unitary tensor functor, then $F$ extends uniquely (up to unitary monoidal equivalence) to a unitary tensor functor $F''\colon\Ccal''\ra\Dcal$. To construct this functor define $F'((U_1,\ldots, U_n)):= F(U_1)\oplus\ldots\oplus F(U_n)$ and on morphisms $F'((T_{ij})_{ij}):=(F(T_{ij}))_{ij}$. If $(U,p)\in\Ob(\Ccal'')$, then $F'(p)$ is a projection in $\End_\Dcal(F(U))$, so there exists $V\in\Ob(\Dcal)$ and an isometry $v\in\Hom_{\Dcal}(V,F(U))$ such that $vv^*=F(p)$ and $v^*v=\iota_V$. Define $F''((U,p)):=V$ and for $p'Tp\in\Hom_{\Ccal''}((U,p),(U',p'))$ let $F''(p'Tp):= v'^*F'(p'Tp)v=v'^*F'(T)v$. The tensor and involutive structure are again defined in the obvious way.\\
Note that in both steps of this extension of $F$ one has to make a choice of objects, a different choice leads to an equivalent functor. Furthermore if both $\Ccal$ and $\Dcal$ are not necessarily closed under direct sums and subobjects and $F\colon\Ccal\ra\Dcal$ is a unitary tensor functor, then $F$ extends to a functor $\colon\Ccal''\ra\Dcal''$. This extension is constructed by applying the universal property to $i\circ F$, where $i:\Dcal\ra\Dcal''$ is the inclusion. The properties ``fully faithful'' and ``essentially surjective'' are preserved under this extension of tensor functors.
\end{Rem} 

\begin{Def}
Let $\Ccal$ be a strict \Cstar tensor category and $U\in\Ob(\Ccal)$. Then $\overline{U}\in\Ob(\Ccal)$ is called {\it conjugate} to $U$ if there exist $R\in\Hom(\unit,\overline{U}\ten U)$ and $\overline{R}\in\Hom(\unit,U\ten\overline{U})$ such that the compositions
\[
\xymatrix{ U\ar[r]^-{\iota\ten R}& U\ten\overline{U}\ten U\ar[r]^-{\overline{R}^*\ten\iota}& U};\qquad 
\xymatrix{\overline{U}\ar[r]^-{\iota\ten\overline{R}}& \overline{U}\ten U\ten\overline{U} \ar[r]^-{R^*\ten\iota} &\overline{U}}
\]
are the identity morphisms. If every object in $\Ccal$ has a conjugate object then $\Ccal$ is {\it rigid}. We say that the pair $(R,\overline{R})$ {\it solves the conjugate equations for} $U$. If $(R,\overline{R})$ is of the form
\[
R=\sum_k (\overline{w_k}\ten w_k)R_k, \qquad \overline{R}=\sum_k (w_k\ten\overline{w_k})\overline{R_k},
\]
where for all $k$ the objects $U_k\in\Ob(\Ccal)$ are simple, $\|R_k\|=\|\overline{R_k}\|$ and $w_k\in\Hom(U_k,U)$ are isometries such that $\sum_k w_kw_k^*=\iota_U$, then $(R,\overline{R})$ is called a {\it standard solution} of the conjugate equations. 
\end{Def}

If an object has a conjugate it also admits a standard solution of the conjugate equations \cite[\textsection 2]{NeshveyevTuset13}. Furthermore, if $(R,\overline{R})$ and $(R',\overline{R}')$ are both standard solutions for $(U,\overline{U})$ respectively $(U,\overline{U}')$, then there exists a unitary $T\in\Hom(\overline{U},\overline{U}')$ such that $R'=(T\ten\iota)R$ and $\overline{R}'=(\iota\ten T)\overline{R}$ \cite[Prop. 2.2.13]{NeshveyevTuset13}. We will only deal with rigid categories. 

\begin{Def}
Suppose that $U\in\Ob(\Ccal)$ and $(R,\overline{R})$ is a standard solution of the conjugate equations for $U$. For $T\in \End(U)$ let $\Tr_U(T)$ be the composition
\[
\xymatrix{\unit\ar[r]^-{R}&\overline{U}\ten U\ar[r]^{\iota\ten T}&\overline{U}\ten U\ar[r]^-{R^*}& \unit}.
\]
This functional $\Tr_U\colon\End(U)\ra\Com$ is called the {\it categorical trace} of $U$. Note that from the remark above it is immediate that the categorical trace is independent on the choice of the standard solution. 
\end{Def}

\begin{Prop}[\protect{\cite[Thm. 2.2.16]{NeshveyevTuset13}}]
Let $U\in\Ob(\Ccal)$, then $\Tr_U\colon\End(U)\ra\Com$ is a tracial, positive and faithful functional. Furthermore $\Tr_U(T)=\overline{R}^*(T\ten\iota)\overline{R}$ for any standard solution $(R,\overline{R})$ of the conjugate equations for $U$.
\end{Prop}

\begin{Def}
Two objects $U,V\in \Ob(\Ccal)$ are {\it isomorphic} if there exists an isomorphism in $\Hom(U,V)$.  We write $[U]$ for the equivalence class of objects isomorphic to $U$. Denote by $K^+[\Ccal]$ the {\it fusion semiring} of $\Ccal$, it is the universal semiring ring generated by the equivalence classes $[U]$ of objects $U\in \Ob(\Ccal)$ with sum and product given by
\[
[U]+[V] := [V\oplus V], \qquad [U][V] := [U\ten V].
\]
Note that there is no need to define a subtraction as we define a semiring.
\end{Def}

Before we define a $SU(d)$-type category let us first say something about the representations of the special unitary group. Details can be found in lots of books, e.g., \cite{FultonHarris}. To avoid trivialities we will always assume that $d\geq2$. We have the fundamental (or defining) representation of $SU(d)$ on $V:=\Com^d$ by letting the group elements act on vectors of $V$ in the straightforward way. By the highest weight classification of irreducible representations of $SU(d)$, we can classify the irreducible representations by the tuples 
\[
\Lambda_d:=\{\lambda=(\lambda_1,\ldots,\lambda_{d-1})\in\Nat^{d-1}\,:\, \lambda_1\geq\lambda_2\geq\ldots\geq\lambda_{d-1}\}.
\]
In this paper we use the convention $\Nat:=\{0,1,2,\ldots\}$. We denote $V_\lambda$ for the irreducible representation corresponding to $\lambda$. For $\lambda\in\Nat^{d-1}$ let $|\lambda|:=\lambda_1+\ldots+\lambda_{d-1}$. It can be shown that any irreducible representation $V_\lambda$ is contained in the tensor product $V^{\ten |\lambda|}$. Another special fact for $SU(d)$ is that the $d$-th anti-symmetric tensor power $\bigwedge^d V$ is isomorphic to the trivial representation. Thus there exists a non-zero map $\Com\ra V^{\ten d}$ intertwining the trivial representation and the $d$-th tensor power of the defining representation.
Given two irreducible representation $V_\lambda$ and $V_\mu$, we can decompose their tensor product representation into irreducible representations. So we have
\[
V_\lambda\ten V_\mu = \bigoplus_{\nu \in\Lambda_d} m_{\lambda,\mu,\nu} V_\nu,
\]
for some multiplicities $m_{\lambda,\mu,\nu}:=\dim \Hom(V_\nu,V_\lambda\ten V_\mu)$. Here $m_{\lambda,\mu,\nu} V_\nu:=V_\nu \oplus \ldots \oplus V_\nu$, with $m_{\lambda,\mu,\nu}$ copies.

\begin{Def}\label{Def_SU(d)-type}
A {\it \Cstar tensor category of $SU(d)$-type}, or simply a {\it category of $SU(d)$-type}, is a rigid \Cstar tensor category $\Ccal$ such that the semirings $K^+[\Ccal]$ and $K^+[\Rep(SU(d))]$ are isomorphic. In particular since simple objects can not be further decomposed, this isomorphism maps simple objects onto simple objects. Therefore we can index the equivalence classes of simple objects of a $SU(d)$-type category by the set $\Lambda_d$. An object $X\in\Ccal$ which corresponds to the fundamental representation $[\Com^d]$ of $SU(d)$ will be called the {\it fundamental object} of $\Ccal$. From now on we fix a $SU(d)$-type category with fundamental object $X$ and for every $\lambda\in\Lambda_d$ we fix a simple object $X_\lambda$ corresponding to $\lambda$.
\end{Def}

\begin{Exam}
The object $X_{\{k\}}:=X_{(k,0,\ldots, 0)}$ corresponds to $S^{k}(V)$, the $k$-th symmetric tensor power of the fundamental representation of $SU(d)$ on $V$. For $1\leq k\leq d-1$, the object $X_{\{1^k\}}:=X_{(1,\ldots,1,0,\ldots, 0)}$ corresponds to $\bigwedge^{k}(V)$, the $k$-th antisymmetric tensor power of the fundamental representation.
\end{Exam}

\begin{Exam}\label{conjugate_X}
The conjugate object $\overline{X}$ is isomorphic to $X_{\{1^{d-1}\}}$. Indeed, by the fusion rules of $SU(d)$ it follows that $X\ten X_{\{1^{d-1}\}} \cong \unit \oplus X_{\{21^{d-2}\}}$. Therefore we obtain that $\Hom(\unit,X\ten X_{\{1^{d-1}\}})\neq\{0\}$, from which the claim follows.
\end{Exam}

\begin{Not}
In a not necessarily strict $SU(d)$-type category denote the objects $X^{\ten 1}:=X$ and $X^{\ten n}:=X\ten X^{\ten n-1}$ for $n\geq 2$. Unwrapping this recursive definition gives $X^{\ten n} = X\ten(X\ten (\cdots (X\ten X)\cdots))$ with $n$ factors of $X$.
\end{Not}

\begin{Lemma}\label{Hom_Xm_Xn}
Let $X$ be a fundamental object of a $SU(d)$-type category $\Ccal$. Then 
\begin{equation}\label{eq_Hom_Xm_Xn1}
X^{\ten n} =\bigoplus_{\lambda\in\Lambda_d} m_{\lambda,n} X_\lambda
\end{equation}
and the multiplicities satisfy $m_{\lambda,n}=0$ if $|\lambda|\not\equiv n \pmod{d}$. In particular if $m\not\equiv n \pmod{d}$, then $\Hom(X^{\ten m},X^{\ten n}) = \{0\}$.
\end{Lemma}
\begin{pf}
For $\Rep(SU(d))$ the identity \eqref{eq_Hom_Xm_Xn1} follows for $X=V$ and $X_\lambda=V_\lambda$ from \cite[Prop. 15.25]{FultonHarris}. Since $\Ccal$ is a $SU(d)$-type category it satisfies the same fusion rules as $SU(d)$. 
\end{pf}

It is possible to obtain a new \Cstar tensor category from an existing one by changing the associativity morphisms. This can be done using twists. We will only define twists in the case of a special type of categories, but twists can be defined in other settings as well, for example for representation categories of compact quantum groups see e.g., \cite{NeshveyevTuset13,NeshveyevYamashita13}.

\begin{Def}\label{twisting_category}
Suppose that $\Ccal$ is a strict \Cstar tensor category and $X$ is an object of $\Ccal$. Let $\rho$ be a $d$-th root of unity and assume that $\Hom_{\Ccal}(X^{\otimes m}, X^{\otimes n})=\{0\}$ if $m\not\equiv n\pmod{d}$. Let $\tilde{\Ccal}$ be the category with objects $\{\unit, X, X^{\ten 2}, \ldots \}$ and morphisms $\Hom_{\tilde{\Ccal}}(X^{\ten m}, X^{\ten n}) := \Hom_{\Ccal}(X^{\otimes m}, X^{\otimes n})$. For $a,b,c\in \Nat=\{0,1,2,...\}$ put $\omega(a,b):=\lfloor\frac{a+b}{d}\rfloor - \lfloor\frac{a}{d}\rfloor - \lfloor\frac{b}{d}\rfloor$. Define the morphisms 
\begin{equation}\label{eq-twist1}
\alpha^\rho_{X^{\ten a}, X^{\ten b}, X^{\ten c}}:= \rho^{\omega(a,b)c}\cdot\alpha_{X^{\ten a}, X^{\ten b}, X^{\ten c}}\colon (X^{\ten a}\ten X^{\ten b})\ten X^{\ten c} \ra X^{\ten a}\ten (X^{\ten b}\ten X^{\ten c}).
\end{equation}
It can be checked (see Lemma \ref{lemma_pentagon} below) that the morphisms $\alpha^{\rho}$ satisfy the pentagon axiom. As $\Hom_{\Ccal}(X^{\otimes m}, X^{\otimes n})=\{0\}$ if $m\not\equiv n\pmod{d}$, we have naturality of $\alpha^{\rho}$. Therefore $\alpha^\rho$ define new associativity morphisms on $\tilde{\Ccal}$. Completing $\tilde{\Ccal}$ with respect to subobjects and direct sums and extending $\alpha^\rho$ to this completion gives new associativity morphisms for the \Cstar tensor category generated by $\tilde{\Ccal}$. We denote this category by $\tilde{\Ccal}^\rho$. If $\Ccal$ is generated by $X$ (that is $\Ccal$ is the direct sum and subobject completion of the full subcategory with objects $\{\unit, X, X^{\ten 2},\ldots\}$), we denote the category we obtain in this way by $\Ccal^{\rho}$.
\end{Def}

These associativity morphisms might seem a bit artificial, but one can prove that the functionals $\rho^{\omega(a,b)c}$ as defined above represent all classes in $H^3(\Int/\Int d,\mathbb{T})$, see e.g., \cite[Prop. A.3]{NeshveyevYamashita13}. This cocycle property is exactly needed to make the pentagonal diagram commutative. Note further that in general twisting does not preserve the existence a braiding. 

\begin{Rem}\label{double_twist}
Note that $\rho^{\omega(a,b)c}\rho'^{\omega(a,b)c} = (\rho\rho')^{\omega(a,b)c}$ for all $a,b$ and $c$. So if $\Ccal$ is generated by $X$, we immediately obtain that $(\Ccal^\rho)^{\rho'} \cong \Ccal^{(\rho\rho')}$.
\end{Rem}

\begin{Lemma}\label{lemma_pentagon}
The morphisms $\alpha^{\rho}$ defined in \eqref{eq-twist1} satisfy the pentagon axiom.
\end{Lemma}
\begin{pf}
Since $\alpha$ are associativity morphisms, commutativity of the diagram
\[
\xymatrix{((X^{\ten a}\ten X^{\ten b})\ten X^{\ten c})\ten X^{\ten e}\ar[d]_{\alpha_{1,2,3}^\rho\ten\iota}\ar[rrd]^{\alpha_{12,3,4}^\rho} & &\\
(X^{\ten a}\ten (X^{\ten b}\ten X^{\ten c}))\ten X^{\ten e}\ar[d]_{\alpha_{1,23,4}^\rho} & & (X^{\ten a}\ten X^{\ten b})\ten (X^{\ten c}\ten X^{\ten e})\ar[d]^{\alpha_{1,2,34}^\rho} \\
X^{\ten a}\ten ((X^{\ten b}\ten X^{\ten c})\ten X^{\ten e})\ar[rr]_{\iota\ten\alpha_{2,3,4}^\rho} & & X^{\ten a}\ten (X^{\ten b}\ten (X^{\ten c}\ten X^{\ten e}))
}
\]
is equivalent to
\[
\rho^{\omega(b,c)e}\rho^{\omega(a,b+c)e}\rho^{\omega(a,b)c} = \rho^{\omega(a,b)(c+e)}\rho^{\omega(a+b,c)e}.
\]
For which in turn it is sufficient to prove that
\[
\omega(b,c)e + \omega(a,b+c)e+\omega(a,b)c - \omega(a,b)(c+e) - \omega(a+b,c)e \equiv 0 \pmod{d}.
\]
One can verify directly that this is in fact an equality and not only a congruency.
\end{pf}

\begin{Lemma}\label{twisted_associativity_morphisms}
Suppose that $\Ccal$ is a strict \Cstar tensor category generated by an object $X$ and $\rho$ is a $d$-th root of unity for some $d\geq2$. Let $\alpha$ and $\alpha^\rho$ be the associativity morphisms in $\Ccal$ respectively in $\Ccal^\rho$. Consider for $m,n\geq1$ the associativity morphism $\alpha_{m,n} \colon X^{\ten m}\ten X^{\ten n} \ra X^{\ten m+n}$ in $\Ccal$, defined by the following inductive relations
\begin{align*}
\alpha_{m,n} := \begin{cases}
\iota_{X\ten X}, &\textrm{if } m=n=1;\\
(\iota\ten\alpha_{m-1,1})\alpha_{X,X^{\ten m-1},X}, & \textrm{if } m\geq 2,\, n=1;\\
\alpha_{m+1,n-1}\circ(\alpha_{m,1}\ten \iota^{\ten n-1})\circ\alpha^{-1}_{X^{\ten m},X,X^{\ten n-1}}, &\textrm{if } m\geq 1,\, n\geq 2.
\end{cases}
\end{align*}
Define similarly the morphisms $\alpha^\rho_{m,n}$ in $\Ccal^\rho$. Then it holds that
\[
\alpha^\rho_{m,n} = \rho^{n\lfloor\frac{m}{d}\rfloor} \alpha_{m,n}.
\]
\end{Lemma}
\begin{pf}
Let us prove this lemma by induction on $m$ and $n$. If $m=n=1$, the lemma is trivial. Suppose that $n=1$. Note that because $d\geq 2$ it holds that $\lfloor\frac{1}{d}\rfloor =0$. So by definition of the twist we obtain
\[
\alpha^{\rho}_{X,X^{\ten m-1},X} = \rho^{\lfloor\frac{m}{d}\rfloor - \lfloor\frac{m-1}{d}\rfloor} \alpha_{X,X^{\ten m-1},X}
\]
as a map $(X^{\ten m}) \ten X\ra X\ten((X^{\ten m-1})\ten X)$. Proceeding by induction on $m$ it follows that
\begin{align*}
\alpha_{m,1}^\rho &= (\iota\ten\alpha_{m-1,1}^\rho) \alpha_{X,X^{\ten m-1},X}^\rho\\
&= \rho^{\lfloor\frac{m-1}{d}\rfloor}\rho^{\lfloor\frac{m}{d}\rfloor - \lfloor\frac{m-1}{d}\rfloor}(\iota\ten\alpha_{m-1,1})\alpha_{X,X^{\ten m-1},X}\\
&= \rho^{\lfloor\frac{m}{d}\rfloor}\alpha_{m,1}
\end{align*}
and the lemma is proved for $n=1$. Now suppose that $n>1$. By the definition and induction hypothesis, it holds that s
\begin{align*}
\alpha^{\rho}_{m,n} &= (\rho^{(n-1)\lfloor\frac{m+1}{d}\rfloor}\alpha_{m+1,n-1}) (\rho^{\lfloor\frac{m}{d}\rfloor}\alpha_{m,1}\ten\iota^{\ten n-1}) (\rho^{-(\lfloor\frac{m+1}{d}\rfloor - \lfloor\frac{m}{d}\rfloor)(n-1)}\alpha^{-1}_{X^{\ten m}, X, X^{\ten n-1}})\\
&= \rho^{n\lfloor\frac{m}{d}\rfloor}\alpha_{m,n},
\end{align*}
as desired.
\end{pf}

\section{Hecke algebras}\label{Hecke_algebras}
In this section we will briefly recall some results about Hecke algebras which will be used later when considering $SU(d)$-type categories. More about Hecke algebras can be found in e.g., \cite{Wenzl88}.

\begin{Def}\label{Def_Hecke_algebra}
Given $n\in\Nat$ and $q\in\Com$, define the Hecke algebra $H_n(q)$ to be the unital algebra generated by the $n-1$ elements $g_1,\ldots,g_{n-1}$ which satisfy the following three relations
\begin{align}
&g_ig_j=g_jg_i & &\textrm{if } |i-j|\geq 2; \label{Eq_Def_Hecke1}\\
&g_ig_{i+1}g_i = g_{i+1}g_ig_{i+1}& &\textrm{for } i=1,\ldots,n-2; \label{Eq_Def_Hecke2}\\
&g_i^2 = (q-1)g_i + q& &\textrm{for } i=1,\ldots,n-1. \label{Eq_Def_Hecke3}
\end{align}
\end{Def}

Note that if $q\neq 0$ we have 
\[
g_i\Big(\frac{1-q}{q} +\frac{1}{q}\,g_i\Big) = \frac{1-q}{q}\,g_i + \frac{1}{q}\, ((q-1)g_i + q) = 1.
\]
So for $q\neq 0$ the elements $g_i$ have inverses. We will denote these by $g_i^{-1}:=\frac{1-q}{q} +\frac{1}{q}g_i$.
Observe that if $q=1$ relation \eqref{Eq_Def_Hecke3} reads as $g_i^2=1$ hence $H_n(1)=\Com[S_n]$, the group algebra of the symmetric group on $n$ elements. So for $q=1$ we obtain a map $S_n \ra H_{n}(q)$, but also for general $q\in\Com$ we can define such a map. 

\begin{Def}\label{map_Sn_in_Hn}
An elementary transposition of $S_n$ is an element of the form $\sigma_i:=(i,i+1)$. Any element $\pi\in S_n$ can be written as a product of elementary transpositions $\pi = \sigma_{i_1}\cdots\sigma_{i_k}$. For a permutation $\pi$ choose such a product of shortest length. The corresponding $k$ will be referred to as the {\it length} of $\pi$, we put $l(\pi):= k$. A product of shortest length will be referred to as a {\it reduced expression} for $\pi$. If $e$ is the identity element of $S_n$ we put $g_e:=1\in H_n(q)$. If $\pi\in S_n$ and $\pi\neq e$, we define $g_\pi:= g_{i_1}\ldots g_{i_k}\in H_n(q)$. From the the lemma below it follows that the element $g_\pi$ is well-defined. 
\end{Def}

\begin{Lemma}[\protect{\cite[\textsection 1.1]{Garsia02}}]\label{basis_Sn}
Let $\pi\in S_n$, define $d_\pi(i):=\#\{1\leq j<i\,:\, \pi(j)>\pi(i)\}$. Then $l(\pi)=\sum_{i=1}^n d_\pi(i)$. Put \[
C_{i,j}:=\begin{cases} 1, &\textrm{if } i\geq j;\\
\sigma_i\cdots\sigma_{j-1}, &\textrm{if } i<j,
\end{cases}
\]
Then $C_{n-d_\pi(n),n}\cdots C_{3-d_\pi(3),3}C_{2-d_\pi(2),2}$ is a reduced expression for $\pi$. Any two reduced expressions for $\pi$ can be transformed in one another by only using the transformations
\begin{align*}
&\sigma_i\sigma_j = \sigma_j\sigma_i& &\textrm{if } |i-j|\geq1;\\ 
&\sigma_i\sigma_{i+1}\sigma_i = \sigma_{i+1}\sigma_i\sigma_{i+1}& &\textrm{for } i=1,\ldots,n-2.
\end{align*}
\end{Lemma}

We can embed $H_n(q)$ into $H_{n+1}(q)$ via the homomorphism $H_n(q)\ni g_i\mapsto g_i\in H_{n+1}(q)$. Iterating this procedure we obtain embeddings $i_{m,n}\colon H_m(q)\ra H_n(q)$ for $m\leq n$. The inductive limit of $(H_{n}(q), i_{m,n})$ is denoted by $H_\infty(q)$. Similarly the {\it shift map} $\Sigma\colon H_n(q)\ra H_{n+1}(q)$, $g_i\mapsto g_{i+1}$ yields another embedding. Unless stated otherwise we will use the first embedding. 
Note that \eqref{Eq_Def_Hecke3} can be rewritten as $(g_i+1)(g_i-q)=0$. So $g_i$ has exactly two spectral values: $-1$ and $q$. In accordance with Kazhdan and Wenzl we define the idempotents\footnote{Note that the other choice of idempotents $e'_i:= \frac{1+g_i}{q+1}$ is also used in the literature.} $e_i:= \frac{q-g_i}{q+1}$. 

\begin{Lemma}\label{Hecke_braiding}
Let $\sigma_{m,n}\in S_{m+n}$ be the permutation defined by
\[
\sigma_{m,n}(i):=\begin{cases} i+n, &\textrm{if } 1\leq i\leq m;\\
i-m, &\textrm{if } m+1 \leq i\leq m+n.\end{cases}
\]
Then 
\begin{equation}\label{eq_Hecke_braiding1}
g_{\sigma_{m,n}}g_i = \begin{cases} g_{i+n}g_{\sigma_{m,n}}, &\textrm{if } 1\leq i\leq m-1;\\
g_{i-m}g_{\sigma_{m,n}}, &\textrm{if } m+1 \leq i\leq m+n.\end{cases}
\end{equation}
Explicitly 
\begin{align*}
g_{\sigma_{m,n}}&=(g_n g_{n-1}\cdots g_1)(g_{n+1}g_n\cdots g_2)\cdots(g_{n+m-1}g_{n+m-2}\cdots g_m)\\
&= (g_n g_{n+1} \cdots g_{n+m-1})(g_{n-1}g_n\cdots g_{n+m-2})\cdots(g_1 g_2\cdots g_m).
\end{align*}
\end{Lemma}
\begin{pf}
The explicit formulas for $g_{\sigma_{m,n}}$ follow from the reduced expression of $\sigma_{m,n}$ as stated in Lemma \ref{basis_Sn}. To prove \eqref{eq_Hecke_braiding1}, suppose that $m=1$ and $i>1$, then
\[
(g_n\cdots g_1) g_i = g_n\cdots (g_i g_{i-1}g_i) g_{i-2}\cdots g_1 = g_n\cdots (g_{i-1} g_{i}g_{i-1}) g_{i-2}\cdots g_1 = g_{i-1} (g_n\cdots g_1).
\]
Now for $m>1$ and $i>m$ the statement follows from the case $m=1$, induction on $m$ and the explicit formula of $g_{\sigma_{m,n}}$. For $i<m$ we use the other expression of $g_{\sigma_{m,n}}$. Assume $n=1$, we obtain
\[
(g_1\cdots g_m)g_i = g_1\cdots (g_ig_{i-1}g_i)g_{i-2} \cdots g_m= g_1\cdots (g_{i-1}g_{i}g_{i-1})g_{i-2}\cdots g_m = g_{i-1} (g_1\cdots g_m).
\]
Again the general case follows from this case, induction on $n$ and the explicit formula of $g_{\sigma_{m,n}}$.
\end{pf}

\begin{Not}
For $q\neq 1$ put $[n]_q:= \frac{1-q^n}{1-q} = (1+q+\ldots +q^{n-1})$ and $[n]_1:=n$, it is called the {\it q-analog}, {\it q-bracket} or {\it q-number}. Define the {\it q-factorial}
\[
[1]_q!:=1,\qquad [n]_q!:=[n]_q[n-1]_q!.
\]
\end{Not}

\begin{Def}
Suppose that $q>0$ or $|q|=1$. Define an involution on $H_n(q)$ by $e_i^*:=e_i$ and by antilinear extension. This involution will be called the {\it standard involution} of $H_n(q)$. From now on we will assume that for these values of $q$ the Hecke algebra $H_n(q)$ is equipped with this standard involution. Note that the idempotents $e_i$ become the spectral projections corresponding to the spectral value $-1$ of $g_i$. Furthermore if $q>0$ the elements $g_i$ become self-adjoint.
\end{Def}

\begin{Lemma}\label{properties_An}
Denote $A_n:= \sum_{\sigma\in S_n} g_\sigma$. If $q^m\neq 1$ for $m=1,\ldots, n$ let $E_n:= ([n]_q!)^{-1} A_n$. With this notation the following holds:
\begin{enumerate}[label=(\roman*)]
\item $A_n= (1+g_{n-1}+g_{n-2}g_{n-1}+\ldots+g_1\ldots g_{n-1})A_{n-1}\newline
=A_{n-1} (1+g_{n-1}+g_{n-1}g_{n-2}+\ldots+g_{n-1}\cdots g_{1})\newline 
= (1 +g_{1} + g_{2}g_{1}+\ldots +g_{n-1}\cdots g_{1})\Sigma(A_{n-1})\newline
= \Sigma(A_{n-1})(1 +g_{1} + g_{1}g_{2}+\ldots +g_1\cdots g_{n-1})$;
\item $A_ng_i = g_i A_n = qA_n$, for $i=1,\ldots, n-1$;
\item $E_n$ is a minimal idempotent in $H_n(q)$. If $q\in\Rea$ or $|q|=1$ and $q^m\neq 1$ for $m=1,\ldots n$, it is a projection.
\end{enumerate}
\end{Lemma}
\begin{pf}
(i) From Lemma \ref{basis_Sn} it follows that every element in $\pi\in S_n$ can uniquely be written as $\pi=\sigma_j\sigma_{j+1}\cdots\sigma_{n-1}\pi'$ for some $j\in\{1,\ldots,n\}$ and $\pi'\in S_{n-1}$. For a permutation $\sigma$ denote $\sigma S_n:=\{\sigma\sigma'\,:\, \sigma'\in S_n\}$. By Lemma \ref{basis_Sn}, 
\begin{equation}\label{eq_properties_An1}
S_{n+1} = S_n \cup \sigma_n S_n \cup\ldots\cup (\sigma_1\cdots\sigma_n) S_n.
\end{equation}
Hence by induction the first equality follows. Also $\pi^{-1} = \pi'^{-1}\sigma_{n-1}\cdots\sigma_j$, therefore 
\[
\sum_{\pi\in S_n} \pi = \sum_{\pi\in S_n} \pi^{-1}= \sum_{\pi'\in S_{n-1}} \pi'\Big(\sum_{j=1}^{n} \sigma_{n-1}\cdots\sigma_{j+1}\sigma_j\Big).
\]
Now by induction the second equality in (i) follows. The third and fourth equalities can be proved similarly (one can use the map $\sigma_i\mapsto \sigma_{n-i}$).\\
Assertions (ii) and (iii) can be found in \cite[\textsection 2]{Pinzari07}. But (ii) can also quickly be derived from (i) and induction and statement (iii) follows again from (ii). 
\end{pf}

\begin{Not}\label{Hecke_morphisms}
Define the maps $\alpha$ and $\beta$ on the generators by 
\begin{align*}
\alpha &\colon H_n(q)\ra H_n(q), & g_i&\mapsto q-1-g_i;\\
\beta &\colon H_n(q)\ra H_n(q^{-1}), & g_i&\mapsto -q^{-1}g_i.
\end{align*}
A simple computation shows that $\alpha$ and $\beta$ respect the defining relations of the Hecke algebras (cf. \eqref{Eq_Def_Hecke1} - \eqref{Eq_Def_Hecke3}) and thus that $\alpha$ and $\beta$ are Hecke algebra morphisms. Furthermore $\alpha\circ\alpha = id$ and $\beta\circ\beta = id$. Note also that $\alpha(e_i)=e_i'=\frac{1+g_i}{q+1}$, where $e_i'$ is the other choice of idempotents, as discussed in \textsection \ref{Hecke_algebras}. 
\end{Not}

\begin{Lemma}\label{properties_Bn}
Suppose that $q\neq 0$. Denote $B_n:= \sum_{\sigma\in S_n} (-q)^{-l(\sigma)} g_\sigma$. If $q^m\neq 1$ for $m=1,\ldots, n$, let $F_n:= ([n]_{\frac{1}{q}}!)^{-1} B_n$. With this notation the following holds:
\begin{enumerate}[label=(\roman*)]
\item $B_n= (1-q^{-1}g_{n-1}+q^{-2}g_{n-2}g_{n-1}+\ldots+(-q)^{-(n-1)}g_{1}\cdots g_{n-1})B_{n-1}\newline
= B_{n-1} (1-q^{-1}g_{n-1}+q^{-2}g_{n-1}g_{n-2}+\ldots+(-q)^{-(n-1)}g_{n-1}\cdots g_{1})\newline 
= (1-q^{-1}g_{1} + q^{-2}g_{2}g_{1}+\ldots+ (-q)^{-(n-1)}g_{n-1}\cdots g_{1})\Sigma(B_{n-1})\newline
= \Sigma(B_{n-1})(1-q^{-1}g_{1}+q^{-2}g_{1}g_{2}+\ldots+(-q)^{-(n-1)}g_{1}\ldots g_{n-1})\newline
= q^{-(n-1)} (1 -qg_{n-1}^{-1} + \ldots + (-q)^{n-1} g_1^{-1}\cdots g_{n-1}^{-1})B_{n-1};$
\item $B_ng_i = g_i B_n = -B_n$, for $i=1,\ldots, n-1$;
\item $F_n$ is a minimal idempotent in $H_n(q)$. If $q\in\Rea\setminus\{0\}$, or $|q|=1$ and $q^m\neq 1$ for $m=1,\ldots n$ it is a projection;
\item $\alpha(B_n)=\mu^{-n(n-1)}A_n$.
\end{enumerate}
\end{Lemma}
\begin{pf}
It is immediate that $\beta(A_n)=B_n$. Thus all assertions except from the last equality in item (i) and item (iv) follow from the previous lemma. To prove (i), first note that if $i=1,\ldots,n-2$, then $g_i^{-1}B_{n-1} = (\frac{1-q}{q}+\frac{1}{q}g_i)B_{n-1} = -B_{n-1}$. Therefore
\begin{align*}
g_{i}^{-1}\cdots &g_{n-1}^{-1}B_{n-1} = \frac{1}{q}\,g_{i}^{-1}\cdots g_{n-2}^{-1}g_{n-1}B_{n-1} + \frac{1-q}{q}\, g_{i}^{-1}\cdots g_{n-2}^{-1}B_{n-1} \\
&= \frac{1}{q}\,g_{i}^{-1}\cdots g_{n-2}^{-1} g_{n-1}B_{n-1} + \frac{1-q}{q}\, (-1)^{n-1-i} B_{n-1}\\
&= \Big(q^{-(n-i)}g_i\cdots g_{n-1} + \frac{1-q}{q}\, q^{-(n-i-1)} g_{i+1}\cdots g_{n-1} + \frac{1-q}{q}\, q^{-(n-i-2)} (-1) g_{i+2}\cdots g_{n-1}\\ 
&\qquad + \frac{1-q}{q}\, q^{-(n-i-3)} (-1)^2 g_{i+3}\cdots g_{n-1} +\ldots +  \frac{1-q}{q}\,(-1)^{n-1-i} \Big)B_{n-1}
\end{align*}
Hence
\begin{align}
&\big(g_1^{-1}\cdots g_{n-1}^{-1} + (-q)^{-1}g_2^{-1}\cdots g_{n-1}^{-1} +\ldots + (-q)^{-(n-2)}g_{n-1}^{-1} + (-q)^{-(n-1)}\big)B_{n-1} = \label{eq_properties_Bn1}\\
&\Big(q^{-(n-1)}g_1\cdots g_{n-1} + \frac{1-q}{q}\, \big( q^{-(n-2)} g_{2}\cdots g_{n-1} - q^{-(n-3)} g_{3}\cdots g_{n-1} +\ldots + (-1)^{n-2}\big)\notag \\
&+ (-q)^{-1} \big(q^{-(n-2)}g_2\cdots g_{n-1} + \frac{1-q}{q}\, \big( q^{-(n-3)} g_{3}\cdots g_{n-1} - q^{-(n-4)} g_{4}\cdots g_{n-1} +\ldots + (-1)^{n-3}\big)\big)\notag \\
& +  \ldots + (-q)^{-(n-2)}\big(q^{-1}g_{n-1} + \frac{1-q}{q}\big) + (-q)^{-(n-1)}\Big)B_{n-1}. \notag
\end{align}
Gathering all terms $g_i\cdots g_{n-1}$, the constant in front of $g_i\cdots g_{n-1}$ becomes
\begin{align*}
&\frac{1-q}{q}\,\Big( (-1)^i q^{-(n-i)} + (-q)^{-1}(-1)^{i-1}q^{-(n-i)} + \ldots + (-q)^{-(i-2)}q^{-(n-i)}\Big) + (-q)^{-(i-1)} q^{-(n-i)} \\
&= (1-q)(-1)^i \big(q^{-(n+1-i)} + q^{-(n+1-i)-1} +\ldots + q^{-(n-1)}\big) - (-1)^i q^{-(n-1)} \\
&= (-1)^i (-q)q^{-(n+1-i)} + (-1)^{i}q^{-(n-1)} - (-1)^i q^{-(n-1)} \\
&= (-1)^{i+1} q^{-(n-i)}.
\end{align*}
In the second last equality above we use the fact that we have an alternating sum. We thus obtain that \eqref{eq_properties_Bn1} equals
\begin{align*}
&\big((-1)^{1+1}q^{-(n-1)} g_1\cdots g_{n-1} + (-1)^{2+1}q^{-(n-2)}g_2\cdots g_{n-1} + \ldots + (-1)^{n+1}q^{-(n-n)}\big) B_{n-1}\\
&= (-1)^{n-1}\big((-q)^{-(n-1)} g_1\cdots g_{n-1} + (-q)^{-(n-2)}g_2\cdots g_{n-1} + \ldots + 1\big) B_{n-1},
\end{align*} 
which by the first equality of item (i) in this lemma gives the desired result.\\
We prove (iv) by induction. The case $n=2$ is easy, as 
\[
\alpha(B_2)=\alpha(1-q^{-1}g_1) = 1-q^{-1}(q-1-g_1) = q^{-1}(1+g_1)=\mu^{-2(2-1)}A_2.
\]
To prove the induction step, first note that $\alpha(-qg_i^{-1}) = \alpha(-1+q-g_i) = -1+q -q+1+g_i =g_i$. Therefore using part (i) of this lemma, the induction hypothesis and Lemma \ref{properties_An}, we get
\begin{align*}
\alpha(B_{n+1})&= \alpha\big(q^{-n} (1 -qg_n^{-1} + \ldots + (-q)^n g_1^{-1}\cdots g_n^{-1})B_n\big)\\
&= q^{-n} \big(1 + \alpha(-qg_n^{-1}) +\ldots +  \alpha((-q)^n g_1^{-1}\cdots g_n^{-1})\big) \mu^{-n(n-1)}A_n\\
&= \mu^{-(n+1)n}(1 + g_n +\ldots +  g_1\cdots g_n) A_n\\
&= \mu^{-(n+1)n}A_{n+1},
\end{align*}
as desired.
\end{pf}

\begin{Def}
A {\it trace} $\tr$ on the Hecke algebra $H_\infty(q)$ is a linear functional $\tr\colon H_\infty(q)\ra\Com$ such that $\tr(ab)=\tr(ba)$ for all $a,b\in H_\infty(q)$ and $\tr(1)=1$. The trace is called a {\it Markov trace} if there exists an $\eta\in\Com$ such that for all $n\in\Nat$ and $x,y\in H_n(q)\subset H_\infty(q)$ the equality $\tr(xe_ny)=\eta\tr(xy)$ holds. We will refer to this identity as the {\it Markov property}. It is known that for each $\eta\in\Com$ there exists a Markov trace with $\tr(e_1)=\eta$, for a proof of this fact see \cite[Thm. 5.1]{Jones87}. 
\end{Def}

\begin{Lemma}\label{equivalent_Markov_traces}
Let $\tr$ be a Markov trace on $H_\infty(q)$ and $\varphi\colon H_\infty(q)\ra\Com$ be a functional with the Markov property such that $\tr(e_1)=\varphi(e_1)$, then $\tr=\varphi$. In particular $\varphi$ is tracial.
\end{Lemma}
\begin{pf}
Any element $x\in H_n(q)\subset H_\infty(q)$ can be written as $x=x_1+x_2e_{n-1}x_3$ for some $x_1,x_2,x_3\in H_{n-1}(q)$. Now for $\psi=\tr$ and $\psi=\varphi$ it holds
\[
\psi(x) = \psi(x_1)+\psi(x_2e_{n-1}x_3) =  \psi(x_1)+\psi(e_1)\psi(x_2x_3)
\]
and the lemma follows by induction to $n$ and the fact that $H_\infty(q)=\bigcup_n H_n(q)$.
\end{pf}

\section{Computations in $\Rep(SU_\mu(d))$}\label{computations}
In this section we will make some computations in the category $\Rep(SU_\mu(d))$ (for $\mu\in(0,1]$) which will be needed later on. The results are analogous to \cite{Pinzari07}, but in that paper a different representation of the Hecke algebra in $\End_{\Rep(SU_\mu(d))}(\Hcal^{\ten n})$ is used. See Remark \ref{different_representations} for a short discussion on these two different representations.\\
Since the representation category of a $q$-deformed Lie group is very similar to the representation category of the Lie group itself (cf. \cite[\textsection 10.1]{ChariPressley95}), it is immediate that $\Rep(SU_\mu(d))$ is a $SU(d)$-type category.

\begin{Not}\label{Def_rep_SUmud}
Consider the \Cstar tensor category $\textrm{Hilb}_{\textrm{f}}$, with objects all finite dimensional Hilbert spaces and the collection of morphisms between two objects is given by all linear maps between the corresponding Hilbert spaces. Let $\Hcal:=\Com^d\in \Ob(\textrm{Hilb}_{\textrm{f}})$ and let $\{\psi_i\}_{i=1}^d$ be an orthonormal basis in $\Hcal$. Jimbo and Woronowicz defined the following representation of the Hecke algebra $H_n(q)$. Let $q:=\mu^2$. Define the map $T\in\End(\Hcal\ten\Hcal)$ by
\[
T(\psi_i\ten\psi_j):=\begin{cases} (q-1)\psi_i\ten\psi_j + \mu\psi_j\ten \psi_i, &\textrm{if } i<j;\\
q\psi_i\otimes \psi_j, &\textrm{if } i=j;\\
\mu\psi_j\otimes \psi_i, &\textrm{if } i>j.
\end{cases}
\]
Then a straightforward computation shows that 
\[
\eta\colon H_n(q)\ra\End(\Hcal^{\ten n}), \qquad g_i\mapsto \iota^{\ten i-1}\ten T\ten\iota^{n-i-1}
\]
defines a representation of $H_n(q)$. If it is necessary to keep track of $n$ we write $\eta_n$ for this representation. The action of the idempotents $e_i$ corresponds to the linear map
\begin{equation}\label{eq_action_idempotents1}
\frac{q-T}{q+1}\,(\psi_i\ten\psi_j)=\begin{cases} \frac{1}{q+1}\,(\psi_i\ten\psi_j - \mu \psi_j\otimes \psi_i), &\textrm{if } i<j;\\
0, &\textrm{if } i=j;\\
\frac{1}{q+1}\,(q\psi_i\otimes \psi_j - \mu\psi_j\otimes \psi_i), &\textrm{if } i>j.
\end{cases}
\end{equation}
To define the category $\Rep(SU_\mu(d))$ we also need an embedding $\Com\hookrightarrow \Hcal^{\ten d}$, corresponding to the morphism intertwining the trivial representation of $SU(d)$ on $\Com$ with the $d$-th tensor power of the standard representation on $\Hcal^{\ten d}$. Up to a normalization the following element in $\Hcal^{\ten d}$ plays the role of this embedding $\Com\ra\Hcal^d$
\begin{equation}\label{Eq_def_S}
S:= \sum_{\sigma\in S_d} (-\mu)^{-l(\sigma)}\psi_{\sigma(d)}\ten\cdots\ten\psi_{\sigma(1)}. 
\end{equation}
Here $l(\sigma)$ denotes the length of $\sigma$, see Definition \ref{map_Sn_in_Hn}. We write $S$ both for the element defined in \eqref{Eq_def_S} and for the map $\Com\ra\Hcal^{\ten d}$, $c\mapsto cS$. This element $S$ can be considered as the $q$-deformed determinant. \\
The representation category $\Rep(SU_\mu(d))$ can be described as being the smallest \Cstar tensor category in $\textrm{Hilb}_{\textrm{f}}$ which contains the object $\Hcal$ and the morphisms $S\in\Hom(\Com,\Hcal^{\ten d})$ and $T\in\End(\Hcal^{\ten 2})$.
\end{Not}

Let us compute $\|S\|$. As $\{\psi_i\}_{i=1,\ldots, d}$ is a basis for $\Hcal$ for $\sigma,\sigma'\in S_d$ it follows that $\langle \psi_{\sigma(d)}\ten\cdots\ten\psi_{\sigma(1)}, \psi_{\sigma'(d)}\ten\cdots\ten\psi_{\sigma'(1)} \rangle = \delta_{\sigma,\sigma'}$. So $\|S\|^2 = \langle S,S\rangle = \sum_{\sigma\in S_d} (-\mu)^{-2l(\sigma)}$. By induction, \eqref{eq_properties_An1} and the fact $l(\sigma_i\cdots\sigma_n\sigma) = l(\sigma)+n-i+1$ for $\sigma\in S_n$, it follows that 
\[
\sum_{\pi\in S_{n+1}} q^{l(\pi)} = (1+q+\ldots +q^n)\sum_{\pi\in S_n} q^{l(\pi)} = [n+1]_q [n]_q! = [n+1]_q!
\]
and thus $\|S\|=[d]_\frac{1}{q}!$.\\
Recall the labelling of the simple objects as introduced in Definition \ref{Def_SU(d)-type}. The representation $\eta$ acts as follows.

\begin{Lemma}\label{projection_SUd}
For the representation $\eta\colon H_n(q)\ra\End(\Hcal^{\ten n})$ the morphism $\eta(e_1)\in\End(\Hcal^{\ten 2})$ is the projection onto $\Hcal_{\{1^2\}}$. 
\end{Lemma}
\begin{pf}
Using \eqref{eq_action_idempotents1} we obtain for $i<j$ and a constant $a\in\Com$
\begin{align*}
\eta(e_1)(\psi_i\ten\psi_j + a\psi_j\ten\psi_i) &= \frac{1-\mu a}{q+1}\,(\psi_i\ten\psi_j - \mu\psi_j\ten\psi_i);\\
\eta(e_1)(\psi_i\ten\psi_i) &= 0.
\end{align*}
In particular putting $a=-\mu$ respectively $a= \frac{1}{\mu}$, shows that 
\begin{align*}
\eta(e_1)(\psi_i\ten\psi_j -\mu \psi_j\ten\psi_i)&= \psi_i\ten\psi_j -\mu \psi_j\ten\psi_i;\\
\eta(e_1)(\psi_i\ten\psi_j +\frac{1}{\mu} \,\psi_j\ten\psi_i) &= 0,
\end{align*}
which means that $\eta(e_1)$ is the orthogonal projection onto 
\[
U:=\Span\big(\{\psi_i\ten\psi_j -\mu\psi_j\ten\psi_i\,:\, 1\leq i<j\leq d\}\big).
\]
Since $g_1e_1=e_1g_1$, we have $\eta(g_1)U\subset U$. Thus $U$ is a subobject of $\Hcal^{\ten 2}$ in $\Rep(SU_\mu(d)$. Now note that $\Hcal^{\ten 2}=\Hcal_{\{1^2\}}\oplus\Hcal_{\{2\}}$. Recall that $V_\lambda$ was defined to be the irreducible representation of $SU(d)$ corresponding to $\lambda$. By \cite[\textsection 10.1]{ChariPressley95} the dimensions of $\Hcal_{\lambda}$ are the same as the dimensions of $V_\lambda$. Therefore $\dim(\Hcal_{\{1^2\}})=\frac{1}{2}d(d-1)$ and $\dim(\Hcal_{\{2\}})=\frac{1}{2}d(d+1)$. Note that $\dim(U)=\frac{1}{2}d(d-1)$, therefore $U= \Hcal_{\{1^2\}}$.
\end{pf}

\begin{Rem}\label{different_representations}
Recall the Hecke algebra morphism $\alpha$ of Notation \ref{Hecke_morphisms}. It is immediate that $\eta\circ\alpha$ is also a representation of $H_n(q)$ on $\Hcal^{\ten n}$. This is exactly the representation which Pinzari considers in \cite[\textsection 4]{Pinzari07}. Explicitly $\eta\circ\alpha (g_i) = \iota^{\ten i-1}\ten T'\ten\iota^{n-i-1}$, where
\begin{align*}
T'(\psi_i\ten\psi_j):= ((q-1)\iota - T)(\psi_i\ten\psi_j) = \begin{cases} -\mu\psi_j\ten \psi_i, &\textrm{if } i<j;\\
-\psi_i\otimes \psi_j, &\textrm{if } i=j;\\
(q-1)\psi_i\ten\psi_j -\mu\psi_j\otimes \psi_i, &\textrm{if } i>j.
\end{cases}
\end{align*}
\end{Rem}

For later use we prove the following identities in $\Rep(SU_\mu(d))$. 

\begin{Prop}\label{computations_SUmud}
In $\Rep(SU_\mu(d))$ the following relations hold:
\begin{align}
&S=\eta(B_d)(\psi_d\ten\cdots\ten\psi_1);\label{eq_computations_SUmud0}\\
&S^*S=[d]_{\frac{1}{q}}!\,\iota;\label{eq_computations_SUmud1}\\
&SS^*=\eta(B_d);\label{eq_computations_SUmud2}\\
&(S^*\ten\iota)(\iota\ten S) = (-\mu)^{-(d-1)}[d-1]_{\frac{1}{q}}!\;\iota;\label{eq_computations_SUmud3}\\
&(S^*\ten\iota^{\ten d-1})(\iota^{\ten d-1}\ten S) = (-\mu)^{-(d-1)}\eta(B_{d-1});\label{eq_computations_SUmud4}\\
&\eta(g_1\cdots g_d)(S\ten\iota)=\mu^{d+1}(\iota\ten S).\label{eq_computations_SUmud5}
\end{align}
Here $B_n\in H_n(q)$ is as in Lemma \ref{properties_Bn}.
\end{Prop}
\begin{pf}
As stated before, these identities are closely related to the identities proved by Pinzari in \cite[\textsection 5]{Pinzari07}, in fact one can deduce the relations above to the identities of \cite{Pinzari07}. We will do this first and then we will also show how one can compute everything directly.
We denote Pinzari's $q$-deformed determinant by  $\tilde S:=\sum_{\sigma\in S_d} (-\mu)^{l(\sigma)}  \psi_{\sigma(1)}\ten\cdots\ten \psi_{\sigma(d)}$. Let $r\colon S_d\ra S_d$ be defined by $r(\sigma)(i):=\sigma(d+1-i)$. Then by Lemma \ref{basis_Sn}
\begin{align*}
l(r(\sigma)) &=  \#\{(i,j)\,:\, i<j,\, r(\sigma)(i)>r(\sigma)(j)\}\\
&=  \#\{(i,j)\,:\, i<j,\, \sigma(d+1-i)>\sigma(d+1-j)\}\\
&=  \#\{(i,j)\,:\, i<j,\, \sigma(i)<\sigma(j)\}
\end{align*}
and thus  $l(\sigma)+l(r(\sigma))= d(d-1)/2$. Therefore we obtain
\begin{align*}
\tilde S &= \sum_{\sigma\in S_d} (-\mu)^{l(\sigma)}  \psi_{\sigma(1)}\ten\cdots\ten \psi_{\sigma(d)}\\
&= \sum_{\sigma\in S_d} (-\mu)^{d(d-1)/2-l(r(\sigma))}  \psi_{r(\sigma)(d)}\ten\cdots\ten \psi_{r(\sigma)(1)}\\
&= (-\mu)^{d(d-1)/2} \sum_{\sigma\in S_d} (-\mu)^{-l(\sigma)}  \psi_{\sigma(d)}\ten\cdots\ten \psi_{\sigma(1)}\\
&= (-\mu)^{d(d-1)/2} S.
\end{align*}
With this identity and the properties of $\alpha$ (see Notation \ref{Hecke_morphisms}), we can derive equations \eqref{eq_computations_SUmud0} - \eqref{eq_computations_SUmud5} from the results in \cite[\textsection 5]{Pinzari07}. For example using \cite[Lemma 5.1 b)]{Pinzari07} gives
\begin{align*}
\eta(B_d)\psi_d\ten\cdots\ten\psi_1 &= \mu^{-d(d-1)}(\eta\circ\alpha(A_d)) \psi_d\ten\cdots\ten\psi_1 \\
&= \mu^{-d(d-1)}  (-\mu)^{d(d-1)/2} \tilde S\\
&= \mu^{-d(d-1)}  (-\mu)^{d(d-1)/2} (-\mu)^{d(d-1)/2} S =S.
\end{align*}
Or by \cite[Lemma 5.4]{Pinzari07} 
\begin{align*}
(S^*\ten\iota)(\iota\ten S) &= (-\mu)^{-d(d-1)}(\tilde S^*\ten\iota)(\iota\ten \tilde S)\\
&= (-\mu)^{-d(d-1)}\mu^{d-1}[d-1]_q!\iota\\
&= (-\mu)^{d-1}[d-1]_{\frac{1}{q}}!.
\end{align*}
The other identities can be verified in a similar way, the details are left to the reader. To compute everything directly we start with a general identity. Suppose that $1\leq i_1<i_2<\ldots<i_n\leq d$ and $1\leq j\leq n-1$, then 
\begin{equation}\label{eq_computations_SUmud15}
\eta(g_n \cdots g_j)(\psi_{i_n}\ten\cdots\ten \psi_{i_1})= \mu^{n+1-j}(\psi_{i_n}\ten\cdots\ten\psi_{i_{n+2-j}}\ten\psi_{i_{n-j}}\ten\cdots\ten\psi_{i_1}\ten\psi_{i_{n+1-j}}).
\end{equation}
Now suppose that $\theta\in S_n$. From Lemma \ref{basis_Sn} we have the reduced expression $\theta=(\theta^{-1})^{-1} = (C_{c_n,n}\cdots C_{c_3,3}C_{c_2,2})^{-1}$, where $c_i=i-d_{\theta^{-1}(i)}$. This gives in combination with \eqref{eq_computations_SUmud15} and the fact $l(\theta)=l(\theta^{-1})$, that the following identity holds
\begin{equation}\label{eq_computations_SUmud14}
\eta(g_\theta)(\psi_{i_n}\ten\cdots\ten\psi_{i_1}) = \mu^{l(\theta)}(\psi_{i_{\theta^{-1}(n)}}\ten\cdots\ten\psi_{i_{\theta^{-1}(1)}}).
\end{equation}
Suppose again $1\leq i_1<i_2<\ldots<i_n\leq d$, define $S_{i_n,\ldots,i_1}:= \sum_{\sigma\in S_n} (-\mu)^{-l(\sigma)}\psi_{i_{\sigma(n)}}\ten\cdots\ten\psi_{i_{\sigma(1)}}$. By \eqref{eq_computations_SUmud14} and Lemma \ref{properties_Bn} we get
\begin{align}
\eta(B_n)(\psi_{i_{\theta^{-1}(n)}}\ten\cdots\ten\psi_{i_{\theta^{-1}(1)}}) &= \mu^{-l(\theta)}\eta(B_n)\eta(g_\theta)(\psi_{i_n}\ten\cdots\ten\psi_{i_1}) \notag \\
&= (-\mu)^{-l(\theta)}\sum_{\sigma\in S_n} (-q)^{-l(\sigma)}\eta(g_\sigma)(\psi_{i_n}\ten\cdots\ten\psi_{i_1}) \notag \\
&= (-\mu)^{-l(\theta)}\sum_{\sigma\in S_n} (-\mu)^{-l(\sigma)} (\psi_{i_{\sigma^{-1}(n)}}\ten\cdots\ten\psi_{i_{\sigma^{-1}(1)}})\notag \\
&= (-\mu)^{-l(\theta)}\sum_{\sigma\in S_n} (-\mu)^{-l(\sigma)} (\psi_{i_{\sigma(n)}}\ten\cdots\ten\psi_{i_{\sigma(1)}})\notag \\
&= (-\mu)^{-l(\theta)}S_{i_n,\ldots,i_1}.\label{eq_computations_SUmud11} 
\end{align}
Setting $n=d$, $(i_1,\ldots,i_d)=(1,\ldots d)$ and $\theta=id$ gives $S_{i_d,\ldots,i_1}=S$ and proves \eqref{eq_computations_SUmud0}.

\smallskip
Equation \eqref{eq_computations_SUmud1} is immediate from the norm of $S$. 

\smallskip
Instead of proving \eqref{eq_computations_SUmud2}, we will prove a stronger statement which we will use later in the proof of this proposition. Using the notation introduced above, we will show that 
\begin{equation}\label{eq_computations_SUmud13} 
\sum_{d\geq i_n>\ldots>i_1\geq 1} S_{i_n,\ldots,i_1}S_{i_n,\ldots,i_1}^* = \eta(B_n).
\end{equation}
Suppose that $j_1,\ldots,j_n\in\{i_1,\ldots,i_n\}$. Order the tuple $(j_n,\ldots,j_1)$ in decreasing order so we obtain $k_n\geq \ldots\geq k_2\geq k_1$. Then let $p$ be minimal such that $k_p=j_1$. Then 
\[
\eta(g_{n-p+1}\cdots g_{n-1})(\psi_{k_n}\ten\cdots\ten \psi_{k_1}) = \mu^{p-1}\psi_{k_n}\ten \cdots\ten \psi_{k_{p+1}}\ten \psi_{k_{p-1}}\ten\cdots\ten \psi_{k_1}\ten \psi_{j_1}.
\]
Iterating this procedure, it follows that there exists a $\sigma\in S_n$ and $c\in\Rea\setminus\{0\}$ such that $\psi_{j_n}\ten\cdots\ten\psi_{j_1}=c\eta(g_\sigma)(\psi_{k_n}\ten\cdots\ten\psi_{k_1})$. 
Suppose that $j_{l'}=j_{l''}$ for some $l'\neq l''$, then $k_l=k_{l+1}$ for some $l$. We thus have 
\[
\eta(B_n)(\psi_{k_n}\ten\ldots \ten\psi_{k_1}) = -\eta(B_n)\eta(g_{n-l})(\psi_{k_n}\ten\ldots \ten\psi_{k_1}) = -q\eta(B_n)(\psi_{k_n}\ten\ldots \ten\psi_{k_1}),
\]
where the first equality follows by Lemma \ref{properties_Bn} and the second from the action of $\eta(g_{n-l})$ on $(\psi_{i_n}\ten\ldots \ten\psi_{i_1})$. Recall $q>0$, so in particular $q\neq -1$. Therefore $\eta(B_n)(\psi_{k_n}\ten\ldots \ten\psi_{k_1}) =0$. Now 
\[
\eta(B_n)(\psi_{j_n}\ten\ldots \ten\psi_{j_1})= c\eta(B_n)\eta(g_\sigma)(\psi_{k_n}\ten\ldots \ten\psi_{k_1}) = c(-1)^{l(\sigma)}\eta(B_n)(\psi_{k_n}\ten\ldots \ten\psi_{k_1}) =0.
\]
So 
\[
\ker(\eta(B_n))^{\perp}\subset \Span(\{\psi_{\sigma(i_n)}\ten\cdots\ten\psi_{\sigma(i_1)}\,:\, \sigma\in S_n,\; d\geq i_n>\ldots>i_1\geq 1\}).
\]
Note that $S_{i_n,\ldots,i_1}^*(\psi_{j_n}\ten\ldots\ten\psi_{j_1})=0$ if there does not exist a $\sigma\in S_n$ such that $i_k=j_{\sigma(k)}$ for all $k=1,\ldots,n$. Thus also
\[
\ker\Big(\sum_{i_n>\ldots>i_1} S_{i_n,\ldots,i_1}S_{i_n,\ldots,i_1}^*\Big)^{\perp}\subset \Span(\{\psi_{\sigma(i_n)}\ten\cdots\ten\psi_{\sigma(i_1)}\,:\, \sigma\in S_n,\; d\geq i_n>\ldots>i_1\geq 1\}).
\]
Now by the fact that $l(\sigma)=l(\sigma^{-1})$ and \eqref{eq_computations_SUmud11} we conclude 
\[
\sum_{i_n>\ldots>i_1} S_{i_n,\ldots,i_1}S_{i_n,\ldots,i_1}^* (\psi_{j_{\sigma(n)}}\ten\cdots\ten\psi_{j_{\sigma(1)}}) = (-\mu)^{-l(\sigma)}S_{j_n,\ldots,j_1}
= \eta(B_n) (\psi_{j_{\sigma(n)}}\ten\cdots\ten\psi_{j_{\sigma(1)}}).
\]
Since the vectors $S_{i_n,\ldots,i_1}$ and $S_{j_n,\ldots,j_1}$ are orthogonal if  $(i_n,\ldots,i_1) \neq (j_n,\ldots,j_1)$, it follows that $\sum_{i_n>\ldots>i_1} S_{i_n,\ldots,i_1}S_{i_n,\ldots,i_1}^*$ and $\eta(B_n)$ act the same on the space 
\[
\Span(\{\psi_{\sigma(i_n)}\ten\cdots\ten\psi_{\sigma(i_1)}\,:\, \sigma\in S_n,\; d\geq i_n>\ldots>i_1\geq 1\}).
\]
Hence \eqref{eq_computations_SUmud13} holds. The choice $n=d$ and $(i_d,\ldots,i_1)=(d,\ldots,1)$ gives \eqref{eq_computations_SUmud2}.

\smallskip
For the proof of \eqref{eq_computations_SUmud3} and \eqref{eq_computations_SUmud4} we introduce the following tensors
\begin{align*}
S^{(1)}_j&:=\sum_{\sigma\in S_d,\; \sigma(d)=j}(-\mu)^{-l(\sigma)}\psi_{\sigma(d-1)}\ten\cdots\ten\psi_{\sigma(1)};\\
S^{(2)}_j&:=\sum_{\sigma\in S_d,\; \sigma(1)=j}(-\mu)^{-l(\sigma)}\psi_{\sigma(d)}\ten\cdots\ten\psi_{\sigma(2)}.
\end{align*}
Note that it is immediate that 
\[
S=\sum_{j=1}^d \psi_j\ten S^{(1)}_j = \sum_{j=1}^d S^{(2)}_j\ten \psi_j.
\]
For $\sigma\in S_{d-1}$ and $j\leq d$ define $p(\sigma)\in S_d$ by 
\[
p(\sigma)(i):=\begin{cases} j &\textrm{if } i=1;\\ 
\sigma(i-1) &\textrm{if } \sigma(i-1)<j;\\
\sigma(i-1)+1 &\textrm{if } \sigma(i-1)>j.
\end{cases}
\]
Then $l(p(\sigma))=l(\sigma)+j-1$ and $p\colon S_{d-1}\ra\{\theta\in S_d\,:\, \theta(1)=j\}$ is a bijection. For the tuple $(i_{d-1},\ldots,i_1):= (d,\ldots,j+1,j-1,\ldots 1)$ we then obtain that 
\begin{align}
S_{i_{d-1},\ldots,i_1} &= \sum_{\sigma\in S_{d-1}} (-\mu)^{-l(\sigma)} \psi_{i_{\sigma(d-1)}}\ten \cdots\ten\psi_{i_{\sigma(1)}} \notag \\
&= \sum_{\sigma\in S_d,\; \sigma(1)=j} (-\mu)^{-l(\sigma)+j-1} \psi_{\sigma(d)}\ten \cdots\ten\psi_{\sigma(2)} \notag \\
&= (-\mu)^{j-1} S_j^{(2)}.\label{eq_computations_SUmud12}
\end{align}
Furthermore we have that the map
\[
s\colon\{\sigma\in S_d\,:\, \sigma(d)=j\}\ra \{\sigma\in S_d\,:\, \sigma(1)=j\};\qquad s(\sigma)(i):=\begin{cases} j &\textrm{if } i=1;\\
\sigma(i-1) &\textrm{if } i>1, \end{cases}
\]
is a bijection and one easily checks that $l(s(\sigma)) = l(\sigma)-(d+1)+2j$. It follows that 
\begin{equation}\label{eq_computations_SUmud10}
S^{(1)}_j=(-\mu)^{-(d+1)+2j}S^{(2)}_j.
\end{equation}
Since $\Hcal$ is an irreducible object in $\Rep(SU_\mu(d))$, the morphism $(S^*\ten\iota)(\iota\ten S)$ acts as a scalar. Suppose that $\{\psi_i\}_{i=1}^d$ is an orthonormal basis with respect to the inner product $\langle\cdot,\cdot\rangle$ on $\Hcal$. We obtain
\[
\langle \psi_i, (S^*\ten\iota)(\iota\ten S)\psi_j\rangle_{\Hcal} = \Big\langle\sum_{k=1}^d \psi_k\ten S^{(1)}_k\ten\psi_i, \sum_{k=1}^d \psi_j\ten S^{(2)}_k\ten\psi_k \Big\rangle_{\Hcal^{\ten d+1}} = \langle S^{(1)}_j,S^{(2)}_i\rangle_{\Hcal^{\ten d-1}}.
\]
To compute this scalar $(S^*\ten\iota)(\iota\ten S)$ it thus suffices to compute $\langle S^{(1)}_d,S^{(2)}_d\rangle$. For this we have
\begin{align*}
\langle S^{(1)}_d,S^{(2)}_d\rangle &= (-\mu)^{d+1-2d}\langle S^{(1)}_d,S^{(1)}_d\rangle\\
&= (-\mu)^{-(d-1)} \sum_{\sigma,\theta\in S_d,\;\sigma(d)=\theta(d)=d} (-\mu)^{-l(\sigma)-l(\theta)}\langle \psi_{\sigma(d-1)},\psi_{\theta(d-1)}\rangle\cdots \langle\psi_{\sigma(1)},\psi_{\theta(1)}\rangle \\
&= (-\mu)^{-(d-1)} \sum_{\sigma\in S_d,\;\sigma(d)=d} (-\mu)^{-2l(\sigma)}\\
&= (-\mu)^{-(d-1)} \sum_{\sigma\in S_{d-1}} q^{-2l(\sigma)} = (-\mu)^{-(d-1)} [d-1]_{\frac{1}{q}}!,
\end{align*}
which establishes \eqref{eq_computations_SUmud3}.

\smallskip
Suppose that $\xi_i\in\Hcal$, then 
\begin{align*}
(S^*\ten\iota^{\ten d-1})(\iota^{\ten d-1}\ten S)(\xi_1\ten\cdots\ten\xi_{d-1}) &= \sum_{i,j=1}^d (S_j^{(2)*}\ten\psi_j^*\ten\iota^{\ten d-1})(\xi_1\ten\cdots\ten\xi_{d-1}\ten\psi_i\ten S_i^{(1)})\\
&= \sum_{j=1}^d S_j^{(2)*}(\xi_1\ten\cdots\ten\xi_{d-1})\cdot S_j^{(1)}.
\end{align*}
Thus $(S^*\ten\iota^{\ten d-1})(\iota^{\ten d-1}\ten S) = \sum_{j=1}^d S_j^{(1)}S_j^{(2)*}$, which equals $\sum_{j=1}^d (-\mu)^{-(d+1)+2j}S_j^{(2)}S_j^{(2)*}$ by \eqref{eq_computations_SUmud10}. Using \eqref{eq_computations_SUmud12} this can be written as $\sum_{j=1}^d (-\mu)^{-(d-1)}S_{d,\ldots,j+1,j-1,\ldots, 1}S_{d,\ldots,j+1,j-1,\ldots, 1}^*$. Now we invoke \eqref{eq_computations_SUmud13} to obtain
\[
(S^*\ten\iota^{\ten d-1})(\iota^{\ten d-1}\ten S) = (-\mu)^{-(d-1)}\eta(B_{d-1}).
\]
Thus \eqref{eq_computations_SUmud4} holds.

\smallskip
To prove \eqref{eq_computations_SUmud5} we use \eqref{eq_computations_SUmud0} and Lemma \ref{Hecke_braiding}. We obtain the following 
\begin{align*}
\eta(g_1\cdots g_d)(S\ten\psi_i) &= \eta(g_1\cdots g_d)\eta(B_d)(\psi_d\ten\cdots\ten\psi_1\ten\psi_i) \\
&= \eta(\Sigma(B_d))\eta(g_1\cdots g_d)(\psi_d\ten\cdots\ten\psi_1\ten\psi_i) \\
&= \mu^{d-i}q\mu^{i-1}\,\eta(\Sigma(B_d))(\psi_i\ten\psi_d\ten\cdots\ten\psi_1) \\
&= \mu^{d+1}(\psi_i\ten S),
\end{align*}
which concludes the proof of this proposition.
\end{pf}

\begin{Rem}
From relations \eqref{eq_computations_SUmud3} and \eqref{eq_computations_SUmud4} it follows directly that $R:=\mu^{(d-1)/2}([d-1]_{\frac{1}{q}}!)^{-1/2} S$ and $\overline{R}:=(-1)^{d-1}\mu^{(d-1)/2}([d-1]_{\frac{1}{q}}!)^{-1/2} S$ solve the conjugate equations for $\Hcal$ in $\Rep(SU_\mu(d))$.
\end{Rem}

\section{Representations of Hecke algebras}\label{representation}
Suppose that $\Ccal$ is a strict $SU(d)$-type category, with fundamental object $X$. The aim of this section is to show that one can extract a constant $q$ from $\Ccal$ such that there exists a representation of the Hecke algebra $H_n(q)$ into $\End(X^{\ten n})$. This section is closely related to \cite[\textsection 4]{KazhdanWenzl}. Once we established this representation, we will show that this representation essentially only depends on the constant $q$ and not on the other information of the category $\Ccal$. To obtain this result the Markov traces will be used. 

\begin{Not}
Recall that $V=\Com^d$ is the fundamental representation of $SU(d)$, in $\Rep(SU(d))$ the object $V_{\{1^2\}}$ is a subrepresentation of $V^{\ten 2}$. Therefore if $\Ccal$ is a $SU(d)$-type category, there exists exists a projection $a\in\End(X^{\ten 2})$ and a morphism $v\in\Hom(X_{\{1^2\}}, X^{\ten 2})$ such that $v^*v=id_{X_{\{1^2\}}}$ and $vv^*=a$. We say that $a$ is the projection of $X^{\ten 2}$ onto $X_{\{1^2\}}$. Define the elements $a_{k}:= \iota^{\ten k-1}\ten a\in\End(X^{\ten k+1})$. If $k<n$ we also write $a_k$ for the element $\iota^{k-1}\ten a\ten \iota^{\ten n-k-1}\in\End(X^{\ten n})$. Denote by $\Sigma$ the map $\Sigma(a_{i}):=a_{i+1}$.
\end{Not}

\begin{Lemma}\label{computation_gammaC}
Let $a\in\End(X^{\ten 2})$ be the projection onto $X_{\{1^2\}}\subset X^{\ten 2}$. Put $a_1:=a\ten\iota$ and $a_2:=\iota\ten a$. Then there exists a constant $\gamma\in (0,1]$ such that 
\begin{equation}\label{eq_def_gammaC1}
a_1a_2a_1 - \gamma a_1 = a_2a_1a_2 - \gamma a_2.
\end{equation}
\end{Lemma}
\begin{pf}
This is a slightly stronger statement than what it is proved in \cite[Prop. 4.2]{KazhdanWenzl} this is due to the fact that $a$ is a projection and not only an idempotent, we will follow the proof by Kazhdan and Wenzl. By the fusion rules of $SU(d)$ we have 
\begin{equation*}
X^{\ten 3} \cong \begin{cases} 
X_{\{2,1\}}\oplus X_{\{2,1\}} \oplus X_{\{3\}} &\textrm{if } d=2;\\
X_{\{1^3\}} \oplus X_{\{2,1\}} \oplus X_{\{2,1\}} \oplus X_{\{3\}} &\textrm{if } d\geq3.
\end{cases}
\end{equation*}
Therefore
\begin{equation}\label{eq_computation_gammaC2}
\End(X^{\ten 3}) \cong \begin{cases} 
M_2(\Com) \oplus \Com &\textrm{if } d=2;\\
\Com \oplus M_2(\Com) \oplus \Com &\textrm{if } d\geq3.
\end{cases}
\end{equation}
We now only consider the case $d\geq 3$, the case $d=2$ is similar. By the fusion rules of $SU(d)$ it follows that $X_{\{1^3\}}$ is a subobject of $X_{\{1^2\}}\ten X$, so there exists a projection $p\in\End(X_{\{1^2\}}\ten X)$ and a morphism $v\in\Hom(X_{\{1^3\}},X_{\{1^2\}}\ten X)$ such that $v^*v=id_{X_{\{1^3\}}}$ and $vv^*=p$. Similarly there exists $w\in\Hom(X_{\{1^2\}}, X^{\ten 2})$ such that $w^*w=id_{X_{\{1^2\}}}$ and $ww^*=a$. Then 
\[
a_1|_{X_{\{1^3\}}} = v^*(w^*\ten\iota)a_1(w\ten\iota)v= v^*(w^*\ten\iota)(ww^*\ten\iota)(w\ten\iota)v = v^* id_{X_{\{1^2\}}\ten X} v = id_{X_{\{1^3\}}}.
\]
So $a_1$ acts on ${X_{\{1^3\}}}$ as the identity. Similarly $a_2|_{X_{\{1^3\}}} = id|_{X_{\{1^3\}}}$. Using this terminology of subobjects, $X_{\{3\}}$ is not a subobject of $X_{\{1^2\}}\ten X$ and $X\ten X_{\{1^2\}}$, which implies $a_1|_{X_{\{3\}}} = a_2|_{X_{\{3\}}} = 0$. We have $\dim(\Hom(X_{\{2,1\}},X_{\{1^2\}}\ten X)) = \dim(\Hom(X_{\{2,1\}},X\ten X_{\{1^2\}})) =1$, thus in $\End(X_{\{2,1\}}\oplus X_{\{2,1\}})$ the morphisms $a_i$ act as rank $1$ projections. So using the isomorphism \eqref{eq_computation_gammaC2} there exist rank $1$ projections $f_i\in M_2(\Com)$ such that the projection $a_i\in\End(X^{\ten 3})$ corresponds to $(1,f_i,0)\in \Com \oplus M_2(\Com) \oplus \Com$. Since $\ran(f_1f_2f_1)\subset\ran(f_1)$ there exists a $\gamma_1\in\Com$ such that $f_1f_2f_1=\gamma_1f_1$. Similarly there exists $\gamma_2\in\Com$ such that $f_2f_1f_2 = \gamma_2f_2$. Now 
\begin{equation}\label{eq_computation_gammaC3}
\gamma_1 f_1f_2 = f_1(f_2f_1f_2) = (f_1f_2f_1)f_2 = \gamma_2f_1f_2.
\end{equation}
Hence either $\gamma_1=\gamma_2$ or $f_1f_2=0$ in the latter case we can set $\gamma_1=\gamma_2=0$. Put $\gamma_\Ccal:=\gamma_1$. Because $f_i$ are projections and thus positive, it must hold that $\gamma_\Ccal\in[0,1]$. Since $a_i$ corresponds to $(1,f_i,0)$, \eqref{eq_computation_gammaC3} gives \eqref{eq_def_gammaC1}.
It remains to show that $\gamma\neq 0$, this is non-trivial and makes use of certain projections on objects in $X^{\ten d}$, the proof can be found in \cite[Prop. 4.2]{KazhdanWenzl}.
\end{pf}

\begin{Not}
Put $\gamma_\Ccal$ to be the constant obtained from $\Ccal$ as in the previous lemma. Pick $q_\Ccal$ such that $\gamma_{\Ccal}=\frac{q_{\Ccal}}{(1+q_{\Ccal})^2}$, i.e. such that $q_{\Ccal}+q_{\Ccal}^{-1}=\gamma_\Ccal^{-1}-2$. From this it is clear that $q_{\Ccal}$ is uniquely determined up to $q_{\Ccal}\leftrightarrow q_{\Ccal}^{-1}$. Therefore to fix a unique $q_\Ccal$ we select $q_\Ccal\in\{z\in\Com\,:\, |z|\leq 1,\,\im(z)\geq0\}\cup \{z\in\Com\,:\, |z|< 1,\,\im(z)<0\}$. If it is clear which category $\Ccal$ is considered we will omit the subscript $_\Ccal$ in $q_\Ccal$ and $\gamma_{\Ccal}$. 
\end{Not}

\begin{Rem}
At this point it is not clear why $q_\Ccal$ is indeed an invariant of the category. A priori it might be dependent on the choice of $X$. However this constant is indeed independent, we will say more about this issue later (cf. Remark \ref{uniqueness_constants}).
\end{Rem}

\begin{Lemma}\label{constant_q}
For a $SU(d)$-type category  $\Ccal$ we have $q_\Ccal\in(0,1]\cup\{e^{i\alpha}\,:\, 0<\alpha<\frac{2\pi}{3}\}$.
\end{Lemma}
\begin{pf}
The function $(0,1]\ra [2,\infty)$, $q\mapsto q+q^{-1}$ is a bijection, so for $\gamma\in (0,1/4]$ it holds $q\in(0,1]$. If $\gamma\in(\frac{1}{4},1]$, then write $\gamma=\frac{1}{4}\cos^{-2}(\alpha/2)$ for a unique $\alpha\in(0,\frac{2\pi}{3}]$. We have
\[
\gamma = (e^{i\alpha/2} + e^{-i\alpha/2})^{-2} = \frac{e^{i\alpha}}{(1+e^{i\alpha})^{2}},
\]
which implies that $q=e^{i\alpha}$.  
\end{pf}

\begin{Cor}\label{Representation_Hecke}
The map 
\begin{equation}\label{eq_Def_rep1}
H_{n}(q_\Ccal)\ra\End(X^{\ten n}), \qquad e_i\mapsto a_i
\end{equation}
extends to a $*$-representation of the Hecke algebra $H_n(q_\Ccal)$.
\end{Cor}
\begin{pf}
Since $g_i=q-(q+1)e_i$, in the Hecke algebra $H_n(q)$ the relations \eqref{Eq_Def_Hecke1}, \eqref{Eq_Def_Hecke2} and \eqref{Eq_Def_Hecke3} can equivalently be described in terms of the idempotents $e_i$ by 
\begin{align}
&e_ie_j=e_je_i, & &\textrm{if } |i-j|\geq 2; \label{Eq_Def_Hecke_proj1}\\
&e_ie_{i+1}e_i - \frac{q}{(1+q)^2}e_i = e_{i+1}e_ie_{i+1} - \frac{q}{(1+q)^2}e_{i+1}, & & \textrm{for } i=1,\ldots,n-2; \label{Eq_Def_Hecke_proj2}\\
&e_i^2 = e_i,& &\textrm{for } i=1,\ldots,n-1. \label{Eq_Def_Hecke_proj3}
\end{align}
From this characterization, the fact that $a$ is a projection satisfying \eqref{eq_def_gammaC1} and the choice of $q$ it is immediate that the map \eqref{eq_Def_rep1} extends to a representation of $H_n(q)$. Since $e_i$ is self-adjoint in $H_n(q)$ and $a_i$ is self-adjoint in $\End(X^{\ten n})$ the map is $*$-preserving.
\end{pf}

\begin{Lemma}\label{q_root_unity}
If $q_\Ccal=e^{i\alpha}$ for some $0<\alpha<\pi$, then $q_\Ccal$ is a root of unity.
\end{Lemma}
\begin{pf}
We can write $q=e^{2\pi i\beta}$ for some $0<\beta<\frac{1}{2}$. A representation of $H_n(q)$ into a \Cstar algebra is a \Cstar representation if the idempotents $e_i$ are mapped to projections. Such a representation is called trivial if it is a direct sum of representations $\pi_1$ and $\pi_0$ where $\pi_1\colon e_i\mapsto id$ for all $i$ and $\pi_0\colon e_i\mapsto 0$ for all $i$. If $q$ is not a root of unity, then there exists an $m\in\Nat\setminus\{0\}$ such that $m-1<\frac{1}{\beta}<m$. Now \cite[Prop. 2.9]{Wenzl88} implies that there exist no non-trivial \Cstar representations of $H_n(q)$ for $n>((m+1)/2)^2$. However by Corollary \ref{Representation_Hecke} for each $n$ we do have a non-trivial \Cstar representation. Hence $q$ must be a root of unity.
\end{pf}

\begin{Def}\label{Def_representation}
Suppose that $\Ccal$ is a strict $SU(d)$-type category with fundamental object $X$. Consider for $m \leq n$ the map 
\[
i_{m,n}\colon \End(X^{\ten m})\ra \End(X^{\ten n}), \qquad T\mapsto T\otimes \iota^{\ten(n-m)}.
\]
Clearly if $k\leq m\leq n$, then $i_{m,n}i_{k,m}=i_{k,n}$. Thus the algebraic inductive limit of $(\End(X^{\ten n}),i_{m,n})$ exists, denote this limit by $M_\Ccal$. The representations $\theta_n\colon H_n(q_\Ccal)\ra\End(X^{\ten n})$ obtained from Corollary \ref{Representation_Hecke} satisfy $i_{m,n}\circ\theta_m(x)=\theta_n\circ i_{m,n}(x)$ for all $m,n$ and $x\in H_m(q_\Ccal)$. Thus the collection $\{\theta_n\}_n$ extends to a representation of the inductive limits $\theta_\Ccal\colon H_\infty(q_\Ccal)\ra M_\Ccal$.  We denote $\theta_\Ccal(x)=\theta_n(x)=\iota_{m,n}(\theta_m(x))\in\End(X^{\ten n})$ for $x\in H_n(q)\subset H_\infty(q)$. Again we will write just $\theta$ if no confusion is possible.
\end{Def}

\begin{Prop}\label{Markov_categorical_trace}
Let $R\colon\unit\ra \overline{X}\ten X$, $\overline{R}\colon\unit\ra X\ten \overline{X}$ be a standard solution of the conjugate equations. The categorical trace $\Tr_\Ccal$ on $\Ccal$ induces a Markov trace $\tr_{\Ccal}$ on $H_\infty(q)$ via
\begin{equation}\label{eq_Def_Markov_trace1}
\tr_{\Ccal}(x):=\|R\|^{-2n}\Tr_{X^{\ten n}}(\theta_\Ccal(x)), \qquad (x\in H_{n}(q)\subset H_\infty(q)).
\end{equation}
\end{Prop}
\begin{pf}
Recall that if $R,\overline{R}$ and $S,\overline{S}$ are standard solutions for $U$ respectively $V$ then $(\iota\ten R\ten\iota)S$ and $(\iota\ten \overline{S}\ten\iota)\overline{R}$ are a standard solution for $U\ten V$ \cite[Thm 2.2.16]{NeshveyevTuset13}. It follows immediately that $\Tr_{U\ten V} = \Tr_U (\iota\ten\Tr_{V}) = \Tr_V (\Tr_U\ten\iota)$. Suppose that $x\in H_n(q)$ and $n<m$, then 
\[
(\iota^{\ten m-1}\ten\Tr_{X})(\theta(x)) = (\iota^{\ten m-1}\ten \overline{R}^*) (\theta(x)\ten\iota) (\iota^{\ten m-1}\ten \overline{R}) = \|R\|^2 \theta (x)\in \End(X^{\ten m-1}).
\]
We conclude that $\tr_\Ccal$ is independent of $n$.\\
Since $\Tr_{X^{\ten n}}$ is tracial, it only remains to check that $\tr_\Ccal$ has the Markov property. Since $X$ is simple, there exists a scalar $\lambda\in\Com$ such that $(\iota\ten\Tr)(a_1)=\lambda\iota$. Then $(\iota^{\ten n}\ten\Tr)(a_n)=\lambda\iota^{\ten n}$ and also
\[
\tr_\Ccal(e_1) = \|R\|^{-4} \Tr_{X^{\ten 2}}(\theta(e_1)) = \|R\|^{-4} \Tr_X ((\iota\ten\Tr_X)(a_1)) =  \|R\|^{-4} \Tr_X (\lambda\iota) = \|R\|^{-2}\lambda.
\]
Now suppose that $x,y\in H_n(q)$, then
\begin{align*}
\tr_{\Ccal}(xe_ny) &= \|R\|^{-2n-2}\Tr_{X^{\ten n+1}}(\theta(xe_ny)) \\
&= \|R\|^{-2n-2}\Tr_{X^{\ten n}}\circ (\iota^{\ten n}\ten \Tr_X)\big(\theta(x)\theta(e_n)\theta(y)\big) \\
&= \|R\|^{-2n-2}\Tr_{X^{\ten n}}\big(\theta(x)\cdot (\iota^{\ten n-1}\ten ((\iota\ten\Tr_X)(a_1)))\cdot \theta(y)\big) \\
&= \|R\|^{-2n-2}\Tr_{X^{\ten n}}(\lambda\theta(xy)) \\
&= \tr_\Ccal(e_1)\tr_{\Ccal}(xy).
\end{align*}
Hence $\tr_\Ccal$ is a Markov trace.
\end{pf}

With these Markov traces it is possible to show that the representation of the Hecke algebra is independent of the category $\Ccal$ in the following sense. 

\begin{Thm}[Kazhdan--Wenzl]\label{equivalence_of_representations}
If $\Ccal$ is a strict $SU(d)$-type category, then $q_\Ccal \in(0,1]$ and the Markov trace satisfies $\tr_\Ccal(g_1)=\frac{q_\Ccal^d}{[d]_{q_\Ccal}}$. Therefore the kernel of the representation $\theta_{\Ccal}\colon H_n(q_\Ccal)\ra\End(X^{\ten n})$ depends only on $q_\Ccal$. Furthermore $\theta_{\Ccal}(H_n(q_{\Ccal}))=\End(X^{\ten n})$. 
\end{Thm}
\begin{pf}
Since $\|\theta(x^*x)\|=\|\theta(x)\|^2$, it holds that $\theta(x)=0$ if and only if $\theta(x^*x)=0$. Because the categorical trace $\Tr_{X^{\ten n}}$ is faithful, we obtain that 
\[
\ker(\theta\colon H_n(q)\ra\End(X^{\ten n}))=\{x\in H_n(q)\,:\, \tr_\Ccal(x^*x)=0\}.
\]
To characterize the kernel of $\theta$ by Proposition \ref{Markov_categorical_trace} and Lemma \ref{equivalent_Markov_traces} it suffices to show that $\tr_\Ccal(g_1)$ can be computed in terms of $q$. This is non-trivial and has been done by Kazhdan and Wenzl (see \cite[Thm. 4.1]{KazhdanWenzl}), here they also prove surjectivity of $\theta$. The idea of their proof is to decompose $H_n(q)/I_n^\mu\cong \bigoplus_i M_i(\Com)$ as a direct sum of matrix algebras $M_i(\Com)$. Here $I_n^\mu:=\{x\in H_n(q)\,:\, \tr_\mu(xy)=0 \textrm{ for all } y\in H_n(q)\}$ and $\tr_\mu$ is the unique Markov trace such that $\tr_\mu(g_1)=\mu$. The values $\mu=\frac{q^m}{[m]_q}$ for $m\in\Nat$ play a special role (see \cite[Prop. 3.1]{KazhdanWenzl}). These matrix blocks $M_i(\Com)$ are related to Young diagrams and thus to representations of $SU(n)$. This allows to compare the dimensions of $H_n(q)/I_n^\mu$ and $\End(X^{\ten n})$. From these dimensions one can deduce that $\tr(g_1)=\frac{q^d}{[d]_{q}}$ and that $q$ cannot be a non-trivial root of unity. Therefore by Lemmas \ref{constant_q} and \ref{q_root_unity} it follows that $q\in(0,1]$. Furthermore one can show that $I_n^{\tr_\Ccal(g_1)}=\ker(\theta:H_n(q)\ra\End(X^{\ten n}))$ and by using another dimension argument one has $H_n(q)/\ker(\theta)\cong\End(X^{\ten n})$, thus $\theta$ must be surjective.
\end{pf}

\begin{Rem}
Combining the above theorem, Remark \ref{different_representations} and \cite[Prop. 4.1]{Pinzari07} it follows that for $n>d$ the kernel $\ker(\theta_\Ccal\colon H_n(q)\ra\End(X^{\ten n}))$ equals the ideal generated by the element $B_{d+1}\in H_n(q)$.
\end{Rem}

\section{Categories generated by Hecke algebras}\label{Categories_Hecke_algebras}
In this section we will give a number of technical requirements on \Cstar tensor categories which allow us to prove that a category which satisfies these assumptions is in fact unitarily monoidally equivalent to a twist of $\Rep(SU_\mu(d))$.  In the next section we will use this result to show that all $SU(d)$-type categories are equivalent to $\Rep(SU_\mu(d))$. Furthermore we will show that in a special case these categories admit a braiding.

\begin{Assum}\label{assum_cat}
Assume $\Ccal$ is a strict \Cstar tensor category generated by an object $X$ which satisfies the following requirements:
\begin{enumerate}[label=(\roman*)]
\item there exists a constant $q_\Ccal\in(0,1]$ and a projection $a\in\End(X^{\ten 2})$ such that 
\[
(a\ten\iota)(\iota\ten a)(a\ten\iota) - \frac{q_\Ccal}{(1+q_\Ccal)^2}(a\ten\iota) = (\iota\ten a)(a\ten\iota)(\iota\ten a) - \frac{q_\Ccal}{(1+q_\Ccal)^2}(\iota\ten a);
\]
\end{enumerate}
This requirement defines a representation $\theta_n\colon H_n(q_\Ccal)\ra\End(X^{\ten n})$, $e_i\mapsto \iota^{\ten i-1}\ten a\ten \iota^{\ten n-i-1}$.
\begin{enumerate}[resume, label=(\roman*)]
\item $\theta_n\colon H_n(q_\Ccal) \ra \End(X^{\ten n})$ is surjective;
\item $\ker(\theta_n\colon H_n(q_\Ccal)\ra \End(X^{\ten n}))= \ker(\eta_n\colon H_n(q_\Ccal)\ra \End(\Hcal^{\ten n}))$, here $\eta_n$ is as in Notation \ref{Def_rep_SUmud};
\item there exists an integer $d_\Ccal\geq 2$ and a morphism $\nu\in\Hom(\unit, X^{\ten d_\Ccal})$ such that $\nu^*\nu=\iota$ and $\nu\nu^*=\theta(F_{d_\Ccal})$;
\item there exists a $d_\Ccal$-th root of unity $\omega_\Ccal$ such that $\theta(g_{d_\Ccal}\cdots g_1)(\iota\ten\nu)=\omega_\Ccal q_\Ccal^{(d_\Ccal+1)/2} (\nu\ten\iota)$;
\item $\Hom(X^{\ten m}, X^{\ten n})=\{0\}$, if $m\not\equiv n\pmod{d_\Ccal}$.
\end{enumerate}
We let $\mu_\Ccal\in(0,1]$, $\mu_\Ccal:=q_\Ccal^{1/2}$. If it is clear which category is considered, the subscript $_\Ccal$ will be dropped.
\end{Assum}

\begin{Rem}
The results in Proposition \ref{computations_SUmud} show that $\Rep(SU_\mu(d))$ satisfies the conditions (i) and (iii)-(vi) of the assumption above. The fact that the representation $\eta:H_n(q)\ra \End(\Hcal^{\ten n})$ is surjective follows from Theorem \ref{equivalence_of_representations}. So $\Rep(SU_\mu(d))$ satisfies Assumption \ref{assum_cat}.
\end{Rem}

\begin{Not}
If $\Ccal$ is strict, then $\Ccal^{\rho}$ (see Definition \ref{twisting_category}) is in general not strict. We define $\theta_n(g_i)\in\End_{\Ccal^\rho}(X^{\ten n})$ to be the composition
\[
\xymatrix{X^{\ten n}\ar[r]^-{\alpha} & X^{\ten i-1} \ten (X^{\ten 2}\ten X^{\ten n-i-1}) \ar[r]^-{\beta} & X^{\ten i-1}\ten(X^{\ten 2}\ten X^{\ten n-i-1}) \ar[r]^-{\alpha^{-1}}& X^{\ten n}.}
\]
Here $\alpha$ is the appropriate associativity morphism in $\Ccal^\rho$ and $\beta:=\iota^{\ten i-1}\ten(\theta_2(g_1)\ten\iota^{\ten n-i-1})$.
\end{Not}

As shown in the next proposition the constant $\omega$ behaves nicely with respect to twisting the associativity morphisms of a category $\Ccal$. This proposition will be of importance, because in some cases it implies that we can restrict ourselves to the case $\omega=1$.

\begin{Prop}\label{twist_of_sld_category}
Suppose that $\Ccal$ satisfies the requirements of Assumption \ref{assum_cat} and $\rho$ is a root of unity of order $d_\Ccal$, then in $\Ccal^\rho$ the equality $(\nu^*\ten\iota)\theta(g_{d_\Ccal})\cdots\theta(g_1)(\iota\ten\nu)= \rho^{-1}\omega_\Ccal\mu_\Ccal^{d_\Ccal+1}\iota$ holds. In particular if $\tilde{\Ccal}$ is the strictification of $\Ccal^{\omega_\Ccal}$, then in $\tilde{\Ccal}$ it holds that $(\nu^*\ten\iota)\theta(g_{d_\Ccal}\cdots g_1)(\iota\ten\nu)= \mu_\Ccal^{d_\Ccal+1}\iota$.
\end{Prop}
\begin{pf}
Since in general $\Ccal^{\rho}$ is not strict, consider $(\nu^*\ten\iota)\theta(g_d)\cdots\theta(g_1)(\iota\ten\nu)$ which equals the composition
\begin{align}
&\xymatrix{X=X\ten\iota\ar[r]^-{\iota\ten\nu} & X\ten X^{\ten d}\ar[r]^-{\alpha^\rho_1} & X^{\ten 2}\ten X^{\ten d-1}\ar[rr]^-{\theta(g_1)\ten\iota^{\ten d-1}} && X^{\ten 2}\ten X^{\ten d-1} \ar[r]^-{\alpha^\rho_2}&} \notag \\
&\xymatrix{X \ten (X^{\ten 2}\ten X^{\ten d-2}) \ar[rr]^-{\iota\ten(\theta(g_1)\ten\iota^{\ten d-2})} && \cdots \ar[r]^-{\alpha^\rho_{d}}& X^{\ten d-1}\ten X^{\ten 2} \ar[rr]^-{\iota^{\ten d-1}\ten\theta(g_1)}& &X^{\ten d-1}\ten X^{\ten 2} \ar[r]^-{\alpha^\rho_{d+1}} &}\notag\\
&\xymatrix{X^{\ten d}\ten X \ar[r]^-{\nu^*\ten\iota}& \unit\ten X =X.}\label{eq_twist_of_sld_category1}
\end{align}
Here $\alpha_i^\rho$ are the associativity morphisms in $\Ccal^\rho$. The composition of these morphisms $\alpha^\rho_{d+1}\circ\cdots\circ\alpha^\rho_2\circ\alpha^\rho_1$ equals the associativity morphism $\alpha^\rho\colon X\ten X^{\ten d}\ra X^{\ten d}\ten X$, which by Lemma  \ref{twisted_associativity_morphisms} acts as multiplication by $\rho^{-1\lfloor\frac{d}{d}\rfloor} = \rho^{-1}$.
In $\Ccal$ the associativity morphisms are trivial. Thus if we replace $\alpha^\rho$ by the associativity morphisms $\alpha$ of $\Ccal$, in $\Ccal$ the composition \eqref{eq_twist_of_sld_category1} equals $\mu_\Ccal^{d+1}\omega_\Ccal\iota$ by requirement (v) of Assumption \ref{assum_cat}. Hence in $\Ccal^{\rho}$ the morphism \eqref{eq_twist_of_sld_category1} acts as $\rho^{-1}\mu_\Ccal^{d+1}\omega_\Ccal\iota$, as desired.
\end{pf}

\begin{Not}
Let $\delta_\Ccal:=(\omega_\Ccal\mu_\Ccal^{d_\Ccal+1})^{-\frac{1}{d_\Ccal}}$. Denote 
\[
T_{m,n}:=\delta_\Ccal^{mn}\theta_\Ccal(g_{\sigma_{m,n}})\in\End(X^{\ten m+n}).
\]
Observe the crucial property $T_{1,d_\Ccal}=(\omega_\Ccal\mu_\Ccal^{d_\Ccal+1})^{-1}\theta(g_{d_\Ccal}\cdots g_1)$, which implies that $T_{1,d_\Ccal}(\iota\ten\nu)=\nu\ten\iota$.
\end{Not}

The following proposition is similar to \cite[Prop. 2.2 a)]{KazhdanWenzl}. 

\begin{Prop}\label{braiding}
Suppose that $\Ccal$ satisfies Assumption \ref{assum_cat} and $\omega_\Ccal=\pm1$, then the collection of morphisms $\{T_{m,n}\}_{m,n\in\Nat}$ defines a braiding on the category $\Ccal$. Explicitly,
\begin{align}
T_{k,m+n} &= (\iota^{\ten m}\ten T_{k,n})(T_{k,m}\ten\iota^{\ten n}); \label{eq_braiding1}\\
T_{k+m,n} &=(T_{k,n}\ten \iota^{\ten m})(\iota^{\ten k}\ten T_{m,n}); \label{eq_braiding2}\\
(\beta\ten\alpha)T_{k,m} &= T_{l,n}(\alpha\ten\beta), \qquad \qquad \textrm{for all }\alpha\in\Hom(X^{\ten k}, X^{\ten l}),\, \beta\in\Hom(X^{\ten m}, X^{\ten n}).\label{eq_braiding9}
\end{align}
\end{Prop}

Note that the case $\omega=-1$ can only occur when $d$ is even, because $\omega$ is a $d$-th root of unity.

\medskip
\begin{pf}
From the explicit formulas in Lemma \ref{Hecke_braiding} we obtain the identities 
\[
\Sigma^m(g_{\sigma_{k,n}})g_{\sigma_{k,m}} = g_{\sigma_{k,m+n}}, \qquad g_{\sigma_{k,n}}\Sigma^k(g_{\sigma_{m,n}})=g_{\sigma_{k+m,n}},
\]
from which \eqref{eq_braiding1} and \eqref{eq_braiding2} immediately follow. Denote the morphism $\nu_{m,n}:=\iota^{\ten m}\ten\nu\ten\iota^{\ten n}\in\Hom(X^{\ten m+n}, X^{\ten m+d+n})$. The collection $\{T_{m,n}\}_{m,n}$ satisfies the following relations
\begin{align}
&(T_{m,d}\ten\iota^{\ten n})\nu_{m,n}=\nu_{0,m+n};\label{eq_braiding3}\\
&(T_{d,m}\ten\iota^{\ten n})\nu_{0,m+n}=\nu_{m,n}; \label{eq_braiding4}\\
&\nu_{m,n}^*(T_{d,m}\ten\iota^{\ten n}) = \nu^*_{0,m+n}; \label{eq_braiding5}\\
&\nu_{0,m+n}^*(T_{m,d}\ten\iota^{\ten n}) = \nu^*_{m,n}. \label{eq_braiding6}
\end{align}
Indeed, the case $m=1$ of \eqref{eq_braiding3} follows immediately from
\begin{equation}\label{eq_braiding7}
(\nu^*\ten\iota)\eta(g_d\cdots g_1)(\iota\ten\nu)=\omega\mu^{d+1}\iota
\end{equation}
and the definition of $T_{1,d}$. The case $m>1$ can be proved using induction and \eqref{eq_braiding2}. For \eqref{eq_braiding4} observe that taking the adjoint of \eqref{eq_braiding7} gives
\[
(\iota\ten\nu^*)\eta(g_1\cdots g_d)(\nu\ten\iota)\omega\mu^{d+1}\iota
\]
Here it is crucial that $\omega= \pm1$, otherwise we would have the factor $\overline{\omega}$. From this equation \eqref{eq_braiding4} follows for $m=1$ and the general case can again be proved by induction. The identities \eqref{eq_braiding5} and \eqref{eq_braiding6} follow from respectively \eqref{eq_braiding3} and \eqref{eq_braiding4} by taking conjugates. Again the requirement $\omega=\pm1$ is implicitly used. 
\\
By assumption on $\Ccal$ the map $\theta\colon H_n(q)\ra\End(X^{\ten n})$ is surjective. Combination with Lemma \ref{Hecke_braiding} gives immediately that for all $\alpha\in\End(X^{\ten m})$ and $\beta\in\End(X^{\ten n})$
\begin{equation}\label{eq_braiding8}
T_{m,n}(\alpha\ten\beta)=(\beta\ten\alpha)T_{m,n}.
\end{equation}
It remains to show that \eqref{eq_braiding9} also holds for $\alpha\in\Hom(X^{\ten k},X^{\ten l})$ and $\beta\in\Hom(X^{\ten m},X^{\ten n})$. We may assume that $k=l+pd$ and $m=n+qd$ for some $p,q\in\Int$. We will proceed by induction on $p$ and $q$. The basis case $p=q=0$ is exactly \eqref{eq_braiding8}. So first suppose that $p\geq1$, $q=0$, $\alpha\in\Hom(X^{\ten k},X^{\ten l})$ and $\beta\in\Hom(X^{\ten m},X^{\ten m})$. Then $(\nu\ten \alpha)\in\Hom(X^{\ten k}, X^{\ten l+d})$. Using the induction hypothesis, \eqref{eq_braiding2} and \eqref{eq_braiding4} we have
\begin{align*}
\nu_{m,l}(\beta\ten\alpha)T_{k,m} &= (\beta\ten\nu\ten\alpha)T_{k,m} = T_{l+d,m} (\nu\ten\alpha\ten\beta)\\
&= (T_{d,m}\ten\iota^{\ten l})(\iota^{\ten d}\ten T_{l,m})(\nu\ten\alpha\ten\beta) = (T_{d,m}\ten\iota^{\ten l})(\nu_{0,l+m})T_{l,m}(\alpha\ten\beta)\\
&=\nu_{m,l}T_{l,m}(\alpha\ten\beta).
\end{align*}
Since the map $\Hom(X^{\ten r},X^{\ten u+v})\ra\Hom(X^{\ten r},X^{\ten u+d+v})$, $\gamma\mapsto \nu_{u,v}\circ\gamma$ is injective, it follows by induction that $(\beta\ten\alpha)T_{k,m} =T_{l,m}(\alpha\ten\beta)$. Now suppose that $p<0$, then $(\nu^*\ten a)\in\Hom(X^{\ten k+d},X^{\ten l})$ with a similar argument as above involving the relations \eqref{eq_braiding1} and \eqref{eq_braiding5} one can show that 
\[
(\beta\ten\alpha)T_{k,m}\nu^*_{0,k+m} = T_{l,n}(\alpha\ten\beta)\nu^*_{0,k+m}.
\]
Injectivity of the map $\Hom(X^{\ten u+v},X^{\ten s})\ra\Hom(X^{\ten u+d+v},X^{\ten s})$, $\gamma\mapsto\gamma\circ\nu^*_{u,v}$ closes the induction on $p$. Induction on $q$ is similar and thus \eqref{eq_braiding9} holds.
\end{pf}

\begin{Thm}\label{characterization}
Suppose that $\Ccal$ satisfies the requirements of Assumptions \ref{assum_cat}. Then $\Ccal$ is unitarily monoidally equivalent to $\Rep(SU_{\mu_\Ccal}(d_\Ccal))^{\overline{\omega_\Ccal}}$. 
\end{Thm}

This theorem uses the ideas of monoidal algebras as described by Kazhdan and Wenzl in \cite[\textsection 2]{KazhdanWenzl}. The proof of this theorem is very similar to the proof of \cite[Proposition 2.2 b)]{KazhdanWenzl} and therefore the computational details will be omitted.

\medskip
\begin{pf}
From Proposition \ref{twist_of_sld_category} and Remark \ref{double_twist} it follows that it suffices to consider the case $\omega_\Ccal=1$. The idea of the proof of this theorem is to extend the isomorphisms $\End(X^{\ten n})\ra \End(\Hcal^{\ten n})$ to $\Hom(X^{\ten k}, X^{\ten l})\ra \Hom(\Hcal^{\ten k}, \Hcal^{\ten l})$ by embedding $\Hom(X^{\ten k}, X^{\ten l})$ into $\End(X^{\ten p})$ for some large $p\in\Nat$ using the maps $\alpha\mapsto \alpha\ten\nu$ and $\alpha\mapsto\alpha\ten\nu^*$. For this, suppose that $k,l,m,n,p\in\Nat$ such that $p=m+kd=n+ld$. We will define some subspaces and maps for $\Ccal$. Note that these constructions can of course also be performed in $\Rep(SU_\mu(d))$. Define the map
\begin{align*}
H_p^{m,n}\colon \Hom(X^{\ten m},X^{\ten n})&\ra\End(X^{\ten p}),\\
\alpha\mapsto (\nu^{\ten l}\ten\iota^{\ten n})\alpha((\nu^*)^{\ten k}\ten\iota^{\ten m}) &= \nu^{\ten l}\ten(\nu^*)^{\ten k}\ten\alpha.
\end{align*}
Then cleary $H_p^{m,n}$ is linear. Define the subspace $\Sigma_p^{m,n}\subset\End(X^{\ten p})$ to be 
\[
\Sigma_p^{m,n}:=\{\beta\in\End(X^{\ten p})\,:\, ((\nu\nu^*)^{\ten l}\ten\iota^{\ten n})\beta = \beta ((\nu\nu^*)^{\ten k}\ten\iota^{\ten m}) = \beta\}.
\]

The proof of the following lemmas is omitted, because it is very similar to \cite[\textsection 2]{KazhdanWenzl}, the only additional requirement one has to check is compatibility of the $*$-structure, but this follows directly from the definitions.

\begin{Lemma}
$H_p^{m,n}$ is an isomorphism of $\Hom(X^{\ten m},X^{\ten n})$ onto $\Sigma_p^{m,n}$. Furthermore for $\alpha\in\Hom(X^{\ten m},X^{\ten n})$ and $\beta\in\Hom(X^{\ten n},X^{\ten r})$ the following identities hold
\begin{equation*}
H_p^{m,n}(\alpha)^* = H_p^{n,m}(\alpha^*), \qquad H_p^{n,r}(\beta)\circ H_p^{m,n}(\alpha)=H_p^{m,r}(\beta\circ\alpha).
\end{equation*}
\end{Lemma}

For each $p$, let $\psi_p\colon\End(X^{\ten p})\ra\End(\Hcal^{\ten p})$ be a $*$-isomorphism making the diagram
\[
\xymatrix{H_p(q)\ar[r]^-{\theta_p}\ar[rd]_{\eta_p} & \End(X^{\ten p})\ar[d]^{\psi_p} \\ & \End(\Hcal^{\ten p})}
\]
commute. Such an isomorphism exists, because by assumption and Theorem \ref{equivalence_of_representations} $\theta_p\colon H_p(q)\ra\End(X^{\ten p})$ and $\eta_p\colon H_p(q)\ra\End(\Hcal^{\ten p})$ are surjective and $\ker(\theta) = \ker(\eta)$. Let us write $\kappa:= \|S\|^{-1} S$, where $S\colon \unit\ra\Hcal^{\ten d}$ is the intertwiner defined in \eqref{Eq_def_S}. Because $\nu\nu^*=\theta(F_d)$ and $\kappa\kappa^* = \|S\|^{-2}SS^* = \eta(F_d)$, we have $\psi_{d}(\nu\nu^*)=\|S\|^{-2}SS^*=\kappa\kappa^*$. Define for $m\equiv n\pmod{d}$ the map $\psi_{m,n}$ which is the composition
\[
\xymatrix{\Hom_{\Ccal}(X^{\ten m},X^{\ten n})\ar[r]^-{H_{p,\Ccal}^{m,n}} &\Sigma_{p,\Ccal}^{m,n} \ar[r]^-{\psi_p} & \Sigma_{p,\Rep(SU_\mu(d))}^{m,n} \ar[rr]^-{(H_{p,\Rep(SU_\mu(d))}^{m,n})^{-1}} &&  \Hom_{\Rep(SU_\mu(d))}(\Hcal^{\ten m},\Hcal^{\ten n}).}
\]

\begin{Lemma}\label{fully_faithful}
The morphisms $\psi_{m,n}$ are well-defined (independent of $p$) isomorphisms of linear spaces and satisfy 
\begin{equation*}
\psi_{m,n}(\alpha^*) = \psi_{n,m}(\alpha)^*, \qquad \psi_{n,r}(\beta)\circ\psi_{m,n}(\alpha) = \psi_{m,r}(\beta\circ\alpha). 
\end{equation*}
\end{Lemma}

With these isomorphisms $(\psi_{m,n})_{m,n}$ we are able to define a unitary tensor functor from $\Ccal$ to $\Rep(SU_\mu(d))$. For this consider the full subcategory $\tilde{\Ccal}$ of $\Ccal$ with objects $\Ob(\tilde{\Ccal} ):=\{X^{\ten n}\,:\, n\in\Nat\}$ and $\Dcal$ the full subcategory of $\Rep(SU_\mu(d))$ with objects $\Ob(\Dcal):=\{\Hcal^{\ten n}\,:\, n\in\Nat\}$. Then the completion of $\tilde{\Ccal}$ and $\Dcal$ with respect to direct sums and subobjects equal respectively $\Ccal$ and $\Rep(SU_\mu(d))$.
Define $\tilde{F}\colon\tilde{\Ccal}\ra\Dcal$ by $X^{\ten n}\mapsto\Hcal^{\ten n}$ on objects, $\tilde{F}(\alpha):= \psi_{m,n}(\alpha)$ for morphisms $\alpha\in\Hom(X^{\ten m}, X^{\ten n})$ and $\tilde{F}_0=id$, $\tilde{F}_2=id$. $\tilde{F}(\alpha)$ is well-defined, because by assumption $m\equiv n\pmod{d}$ if $\alpha\neq 0$.

\begin{Lemma}\label{unitary_tensor_functor}
$\tilde{F}$ is a unitary tensor functor.
\end{Lemma}

Clearly $\tilde{F}$ is essentially surjective. Note that Lemmas \ref{fully_faithful} and \ref{unitary_tensor_functor} imply that $\tilde{F}$ is a fully faithful unitary tensor functor. Taking the completions of $\tilde{\Ccal}$ and $\Dcal$ with respect to direct sums and subobjects gives us the categories $\Ccal$ and $\Rep(SU_\mu(d))$. Under this completion $\tilde{F}$ extends uniquely (up to natural unitary isomorphism) to a unitary tensor functor $F\colon\Ccal\ra\Rep(SU_\mu(d))$. Then $F$ is again fully faithful and essentially surjective. So $F$ is a unitary monoidal equivalence, in other words $\Ccal$ is unitarily monoidally equivalent to $\Rep(SU_\mu(d))$.
\end{pf}

\section{Two characterizations of $SU(d)$-type categories}\label{2characterizations}
The aim of this section is to prove the main results of this paper, namely to characterise all $SU(d)$-type categories and to give a condition when it is possible to embed $\Rep(SU_\mu(d))$ in a given \Cstar tensor category. It will be shown that all $SU(d)$-type categories can be classified by a pair $(q,\omega)$ where $q\in(0,1]$ and $\omega$ is a $d$-th root of unity. The requirement for existence an embedding is given by six identities which basically state that if a category satisfies those requirements, there exist a representation of the Hecke algebra, and the twist and solutions of the conjugate equations can be explicitly computed. The proofs of both theorems consist of showing that in both cases the Assumptions \ref{assum_cat} are satisfied allowing to apply Theorem \ref{characterization}. 

\begin{Def}\label{def_twist}
Let $\Ccal$ be a strict $SU(d)$-type category. Since in $\Rep(SU(d))$ the trivial representation $\Com$ is a subrepresentation of $V^{\ten d}$, there exist a morphism $\nu\colon\unit\hookrightarrow X^{\ten d}$, such that $\nu^*\nu= id_{\unit}$ and $\nu\nu^*\in\End(X^{\ten d})$ is a projection. We define the {\it twist} $\tau_\Ccal$ of $\Ccal$ to be the number by which one multiplies in the following composition\footnote{Note that this twist differs a factor $(-1)^d$ from the twist defined in \cite{KazhdanWenzl}.}
\[
\xymatrix{X=X\ten\unit\ar[r]^-{\iota\ten\nu}& X\ten X^{\ten d} \ar[rr]^{\theta(g_d\cdots g_1)} &&X^{\ten d}\ten X \ar[r]^-{\nu^*\ten\iota} & \unit\ten X = X.}
\]
Note that since $X$ is simple, we obtain a scalar. Also $\tau_\Ccal$ is clearly independent of the choice of $\nu$. Again, a priori it is not clear why $\tau_\Ccal$ is independent of the choice of $X$. Fortunately this is the case as we will show later (cf. Remark \ref{uniqueness_constants}). 
\end{Def}

\begin{Lemma}\label{shift}
The following holds: $\theta(g_d\cdots g_1)(\iota\ten\nu) = \tau_\Ccal (\nu\ten\iota)$.
\end{Lemma}
\begin{pf}
Note that $\nu\nu^*\in\End(X^{\ten d})$. From Theorem \ref{equivalence_of_representations} we obtain that there exists a $x\in H_d(q)$ such that $\nu\nu^*=\theta(x)$. By Lemma \ref{Hecke_braiding} it therefore follows that $\theta(g_d\cdots g_1)(\iota\ten\nu\nu^*) = (\nu\nu^*\ten\iota)\theta(g_d\cdots g_1)$. Since $\nu^*\nu=\iota$ we have
\[
\theta(g_d\cdots g_1)(\iota\ten\nu) = \theta(g_d\cdots g_1)(\iota\ten\nu\nu^*)(\iota\ten\nu) = (\nu\nu^*\ten\iota)\theta(g_d\cdots g_1)(\iota\ten\nu)
\]
and the result follows.
\end{pf}

Observe that identity \eqref{eq_computations_SUmud5} of Proposition \ref{computations_SUmud} implies that the twist of $\Rep(SU_\mu(d))$ equals $\mu^{d+1}$. 

\begin{Not}
Since for a strict $SU(d)$-type category $\Ccal$ the constant $q_\Ccal\in(0,1]$, define $\mu_\Ccal\in(0,1]$ to be the positive square root of $q_\Ccal$. 
\end{Not}

\begin{Lemma}\label{Fn_projection}
Let $\Ccal$ be a strict $SU(d)$-type category. For $n\leq d$, the morphism $\theta(F_n)\in\End(X^{\ten n})$ is the projection corresponding to the inclusion $X^{\ten n} \subset X_{\{1^n\}}$. Here $X_{\{1^d\}}:=\unit$, to express the fact that there exists a non-zero morphism $\nu\colon\unit\ra X^{\ten d}$.
\end{Lemma}
\begin{pf}
We proceed by induction. The case $n=2$ is trivial. Suppose that $2\leq n\leq d-1$ and the result holds for $n$. Let $p_k\in\End(X^{\ten k})$ be the projection corresponding to $X_{\{1^k\}}\subset X^{\ten k}$. To prove the induction step we must show that $p_{n+1}=\theta(F_{n+1})$. By the fusion rules of $SU(d)$ we have $X_{\{1^n\}}\ten X \cong X_{\{1^{n+1}\}}\oplus X_{\{2 1^{n-1}\}}$ and $X\ten X_{\{1^n\}} X \cong X_{\{1^{n+1}\}}\oplus X_{\{2 1^{n-1}\}}$. So either $(p_n\ten\iota)(\iota\ten p_n)=p_{n+1}$ or $(p_n\ten\iota)=(\iota\ten p_n)$. Let us argue by contradiction and assume that the second case holds. Let for $i=0,\ldots, n$, $r_i:=\iota^{\ten i}\ten p_n\ten\iota^{n-i}\in\End(X^{\ten 2n})$. From the assumption it follows that $r_i=r_{i+1}$ for all $i$ and therefore we have $r_0=r_n$. In particular $r_0(1-r_n) = r_0(1-r_0)=0$. On the other hand $r_0(1-r_n)$ cannot be zero, because e.g., the non-zero object $X_{\{1^n\}}\ten X_{\{n\}}$ lies in the range of $r_0(1-r_n)$, which yields a contradiction. We conclude that $(p_n\ten\iota)(\iota\ten p_n)=p_{n+1}$.\\
By Lemma \ref{properties_Bn} we have $F_{n+1}F_n=F_{n+1}\Sigma(F_n)=F_{n+1}$, and thus by the induction hypothesis $\theta(F_{n+1})p_{n+1}= \theta(F_{n+1})(p_n\ten\iota)(\iota\ten p_n) = \theta(F_{n+1} F_n\Sigma(F_n))=\theta(F_{n+1})$. Since $X_{\{n+1\}}$ is simple, $p_{n+1}$ is a minimal projection. By the previous calculation $\theta(F_{n+1})$ is a subprojection of $p_{n+1}$. To show that $\theta(F_{n+1})$ equals $p_{n+1}$ it thus suffices to show that $\theta(F_{n+1})\neq 0$. For this we compute
\begin{align*}
(\iota^{\ten n}\ten \tr_\Ccal)(F_{n+1}) &= [n+1]_{\frac{1}{q}}^{-1} \Big(1+\frac{-1}{q}\frac{q^d}{[d]_q} + \big(\frac{-1}{q}\big)^2\frac{q^d}{[d]_q}g_{n-1} + \ldots + \big(\frac{-1}{q}\big)^n g_1\cdots g_{n-1} \frac{q^d}{[d]_q}\Big) F_n\\
&= [n+1]_{\frac{1}{q}}^{-1} (1-\frac{q^d}{[d]_q} \frac{1}{q} [n]_{\frac{1}{q}}) F_n,
\end{align*}
here we used Lemma \ref{properties_Bn}. Since $\frac{q^n}{[n]_q} = \frac{q}{[n]_{\frac{1}{q}}}$ and $n<d$ it follows that 
\[
1-\frac{q^d}{[q]_d} \frac{1}{q} [n]_{\frac{1}{q}} = 1-\frac{q^d}{[d]_q}\,\frac{[n]_q}{q^n} \neq 0. 
\]
By the induction hypothesis $\tr_\Ccal(F_n)\neq 0$, thus $\tr_\Ccal(F_{n+1})\neq 0$ and hence $\theta(F_{n+1})\neq 0$.
\end{pf}

\begin{Cor}\label{root_twist}
Suppose that $\Ccal$ is a strict $SU(d)$-type category, then there exists a $d$-th root of unity $\omega_\Ccal$ such that $\tau_\Ccal = \omega_\Ccal \mu_\Ccal^{d+1}$. 
\end{Cor}
\begin{pf}
First note that Lemma \ref{Hecke_braiding} implies that $g_{\sigma_{k,d}}=g_{\sigma_{k-1,d}}\Sigma^{k-1}(g_{\sigma_{1,d}})$. Combination with the identity $\theta(g_d\cdots g_1)(\iota\ten\nu)=\tau_\Ccal (\nu\ten\iota)$ gives
\[
\theta(g_{\sigma_{k,d}})(\iota^{\ten k}\ten\nu) = \tau_\Ccal \theta(g_{\sigma_{k-1,d}})(\iota^{\ten k-1}\ten\nu\ten\iota).
\]
By induction we obtain for all $k\in\Nat$ 
\[
\theta(g_{\sigma_{k,d}})(\iota^{\ten k}\ten\nu) = \tau_\Ccal^k(\nu\ten\iota^{\ten k}).
\]
Thus in particular
\[
\theta(g_{\sigma_{d,d}})(\iota^{\ten d}\ten\nu) = \tau_\Ccal^d(\nu\ten\iota^{\ten d}).
\]
Multiplying both sides by $(\nu^*\ten\iota^{\ten d})$ gives $(\nu^*\ten\iota^{\ten d})\theta(g_{\sigma_{d,d}})(\iota^{\ten d}\ten\nu) = \tau_\Ccal^d \iota^{\ten d}$. Combination with the above lemma gives that as a morphism in $\End(X^{\ten 2d})$ we have
\[
(\theta(F_d)\ten\theta(F_d))\theta(g_{\sigma_{d,d}})(\theta(F_d)\ten\theta(F_d)) = \tau_\Ccal^d (\theta(F_d)\ten\theta(F_d)).
\]
Theorem \ref{equivalence_of_representations} shows that the representations $\theta$ and $\eta$ are equivalent. In particular this implies that $\tau_\Ccal^d = \tau_{\Rep(SU_\mu(d))}^d = (\mu^{d+1})^d$, which proves the corollary.
\end{pf}

\begin{Rem}
In \cite[Prop. 5.2]{KazhdanWenzl} it is asserted that $\tau_\Ccal = (-1)^d\omega$ for some $d$-th root of unity $\omega$. This is not true as for example the explicit calculation for $SU_\mu(d)$ shows (cf. \eqref{eq_computations_SUmud5}). The mistake in the proof, is that it is claimed that $\theta(g_{\sigma_{d,d}})$ acts as $(-1)^{d^2}$ on the object $X_{\{1^d\}}\ten X_{\{1^d\}}$. 
\end{Rem}

Now all the technical work has been done to give a classification of $SU(d)$-type categories. 

\begin{Thm}\label{classification}
If $\Ccal$ is a $SU(d)$-type category with fundamental object $X$. Then $(\Rep(SU_{\mu_\Ccal}(d)))^{\overline{\omega_\Ccal}}$ is unitarily monoidally equivalent to $\Ccal$. Furthermore $\Ccal$ admits a braiding if $\omega_\Ccal = \pm1$. 
\end{Thm}
\begin{pf}
By Corollary \ref{Representation_Hecke} we have a representation of the Hecke algebra $H_n(q_\Ccal)\ra\End_\Ccal(X^{\ten n})$. By Theorem \ref{equivalence_of_representations} this representation is surjective and depends only on $q_\Ccal$. As the representation $\eta\colon H_2(q_\Ccal)\ra\End_{\Rep(SU_{\mu_\Ccal}(d))}(\Hcal^{\ten 2})$ satisfies that $\eta(e_1)$ is the projection onto $\Hcal_{\{1^2\}}$ (cf. Lemma \ref{projection_SUd}), we obtain that $q_{\Rep(SU_{\mu_\Ccal}(d))}=q_\Ccal$. Then again by Theorem \ref{equivalence_of_representations} $\ker(\eta)=\ker(\theta)$. Lemmas \ref{Hom_Xm_Xn}, \ref{shift}, \ref{Fn_projection} and Corollary \ref{root_twist} show that the other requirements of Assumption \ref{assum_cat} are satisfied. Now Proposition \ref{braiding} and Theorem \ref{characterization} give the result. 
\end{pf}

\begin{Rem}
It can be shown \cite[Rem. 4.4]{NeshveyevYamashita13} that in general a $SU(d)$-type category is not braided; such a category $\Ccal$ admits a braiding if and only if $\omega_\Ccal=\pm1$. 
\end{Rem}

\begin{Rem}\label{uniqueness_constants}
Now we can also prove why the constants $q_\Ccal$ and $\tau_\Ccal$ are independent of the chosen generator $X$ of the category $\Ccal$. By \cite{McMullen84} all automorphisms of $\Rep(SU(d))$ are in 1-1 correspondence with symmetries of the Dynkin diagram of $SU(d)$. This diagram, consisting of $d-1$ nodes $\{1,2,\ldots, d-1\}$ where the nodes $i$ and $i+1$ are connected by a single edge, has exactly two symmetries, namely the identity and the map given on the nodes by $i\mapsto d-i$. So we only have to show that $q_\Ccal$ and $\tau_\Ccal$ are invariant under this second, non-trivial, map. This map induces an automorphism of $U_{\mu}(SU(d))$, the quantum enveloping Hopf algebra of $SU(d)$, given on the generators by $E_i\mapsto E_{d-i}$, $F_i\mapsto F_{d-i}$, $K^{\pm}_i\mapsto K^{\pm}_{d-i}$. In $\Rep(SU_\mu(d))$ it thus maps every object to a conjugate object. Therefore it is sufficient to show that if we would have chosen $\overline{X}$ instead of $X$ as generating object, the resulting constants $q_\Ccal$ and $\tau_\Ccal$ are the same. This is implicitly proved in \cite[\textsection 4.2]{NeshveyevYamashita13}. The idea is the following, suppose that in $\Ccal$ the associativity morphisms are given by a cocycle $\varphi\in H^3(\Int/d\Int,\Torus)$, thus $\alpha\colon(X^{\ten a}\ten X^{\ten b})\ten X^{\ten c} \ra X^{\ten a}\ten (X^{\ten b}\ten X^{\ten c})$ acts as multiplication by $\varphi(a,b,c)$ (in our case $\varphi$ is of the form $\varphi(a,b,c) = \omega_\Ccal^{(\lfloor\frac{a+b}{d}\rfloor - \lfloor\frac{a}{d}\rfloor - \lfloor\frac{a}{d}\rfloor)c}$). We write $\Rep(SU_\mu(d))^\varphi$ for the category $\Rep(SU_{\mu}(d))$ with these new associativity morphisms. Then $\Ccal\cong \Rep(SU_\mu(d))^\varphi$. The map $X\mapsto \overline{X}$ then corresponds to changing the cocycle $\varphi$ to the new one given by $\psi(a,b,c):=\varphi(-a,-b,-c)$. Then one obtains an isomorphism $\theta\colon \Rep(SU_\mu(d))^\varphi \ra \Ccal\ra\Rep(SU_\mu(d))^\psi$. The question is now whether this isomorphism acts trivially on $H^3(\Int/d\Int,\Torus)$. This is indeed the case, since $\varphi=\partial f$, $\psi=\partial g$, where $f(a,b)=\omega^{-\lfloor\frac{a}{d}\rfloor b}$ and $g(a,b)=\omega^{\lfloor\frac{-a}{d}\rfloor b}$ are maps $f,g\colon\Int\times \Int/d\Int\ra\Torus$. Now a direct computation shows that $fg^{-1}$ factors through $\Int/d\Int\times \Int/d\Int$ and thus $\varphi$ and $\psi$ are equivalent cocycles. Thus $\theta$ acts trivially and $\omega_\Ccal$ and $q_\Ccal$ are invariant under $X\mapsto \overline{X}$.

Another (more straightforward) method of proving that those constants are invariant is by explicitly computing everything. This can be done in the following way. We adopt the notation as in \cite[\textsection 2.2]{NeshveyevTuset13} and denote $F\colon\Ccal\ra\Ccal$ for the contravariant tensor functor 
\[
\Ob{\Ccal}\ra\Ob{\Ccal},\qquad U\mapsto \overline{U}; \qquad \qquad \Hom(U,V)\ra\Hom(\overline{V},\overline{U}),\qquad T\mapsto T^\vee,
\]
where $T^\vee:= (\iota\ten \overline{R}_V^*)(\iota\ten T\ten\iota)(R_u\ten\iota)$. Define $F_2(U,V)\colon \overline{V}\ten\overline{U}\ra \overline{U\ten V}$ by the identity 
\[
(F_2(U,V)\ten\iota\ten\iota)(\iota\ten R_U\ten\iota)R_V=R_{U\ten V}.
\]
Put $a^c:=F_2^*(X,X)a^\vee F_2(X,X)$. Then it can be checked that the $a_i^c$ satisfy the relations of $e_i$ in the Hecke algebra $H_n(q_\Ccal)$ and thus we get a representation $\theta^c\colon H_n(q)\ra\End_\Ccal(\overline{X}^{\ten n})$. Hence $q_\Ccal$ is invariant. Now for $\tau_\Ccal$ we define 
\[
\nu^c:=(F_2^*(X,X)\ten\iota^{\ten d-2})\cdots (F_2^*(X^{\ten d-2},X)\ten\iota) F_2^*(X^{\ten d-1},X) (\nu^*)^\vee\colon\unit\ra\overline{X}^{\ten d}.
\]
Then one can verify that $\nu^c$ plays the role of $\nu$ and 
\[
(\nu^{c*}\ten\iota)\theta^c(g_d\cdots g_1)(\iota\ten\nu^c)=\tau_\Ccal\iota,
\]
whence $\tau_\Ccal$ is invariant under the transformation $X\mapsto \overline{X}$.
\end{Rem}

\begin{Rem}
The above theorem says that all $SU(d)$-type categories can be described by a pair $(q,\omega)$, where $q\in(0,1]$ and $\omega$ is a $d$-th root of unity. Namely we have shown that a $SU(d)$-type category $\Ccal$ is isomorphic to $(\Rep(SU_{\sqrt{q}}(d)))^{\overline\omega}$. Now one might wonder if each pair $(q,\omega)$ of this form can be realised by a compact quantum group. This is indeed the case, see \cite{NeshveyevYamashita13}.
\end{Rem}

Inspired by \cite{Pinzari07} we have the following condition for the existence of an embedding of a twist of $\Rep(SU_\mu(d))$ in a \Cstar tensor category $\Dcal$. We use the notation as introduced in Notation \ref{Def_rep_SUmud}.

\begin{Thm}\label{embedding}
Suppose that $\Dcal$ is a strict \Cstar tensor category such that there exists an object $X\in\Ob(\Dcal)$, morphisms $\nu\in\Hom(\unit,X^{\ten d})$, $a\in\End(X^{\ten 2})$, a constant $\mu\in(0,1]$ and a $d$-th root of unity $\omega$ satisfying the following properties:
\begin{align}
&a=a^*=a^2;\label{eq_embedding1}\\
&(a\ten\iota)(\iota\ten a)(a\ten\iota) - \frac{q}{(1+q)^2}(a\ten\iota) = (\iota\ten a)(a\ten\iota)(\iota\ten a) - \frac{q}{(1+q)^2}(\iota\ten a); \label{eq_embedding2}\\
&\nu^*\nu=\iota;\label{eq_embedding3}\\
&\nu\nu^* = \theta(F_d);\label{eq_embedding4}\\
&(\nu^*\ten\iota)(\iota\ten\nu) = \omega(-\mu)^{-(d-1)}[d]_{\frac{1}{q}}^{-1}\,\iota ;\label{eq_embedding5}\\
&\theta(g_d\cdots g_1)(\iota\ten\nu) = \omega \mu^{d+1}(\nu\ten\iota).\label{eq_embedding7}
\end{align}
Here $q:=\mu^2$ and $\theta\colon H_n(q)\ra\End_{\Dcal}(X^{\ten n})$ is the representation of the Hecke algebra as in Corollary \ref{Representation_Hecke}. Let $\Ccal$ be the sub \Cstar tensor category of $\Dcal$ generated by the object $X$ and morphisms $\nu$ and $a$. Then $\Ccal$ is a $SU(d)$-type category and there exists a unique (up to natural unitary isomorphism) unitary tensor functor $F\colon(\Rep(SU_\mu(d)))^\omega\ra\Dcal$ such that $F(\Hcal)=X$ and $F(S)= ([d]_{q}!)^{1/2} \,\nu$, $F(T)= q-(q+1)a$.
\end{Thm}
\begin{pf}
First note that equations \eqref{eq_embedding1} and \eqref{eq_embedding2} together with Corollary \ref{Representation_Hecke} imply that we have a $*$-representation $\theta\colon H_n(q)\ra\End_\Dcal(X^{\ten n})$. Therefore the identities \eqref{eq_embedding4} and \eqref{eq_embedding7} make sense. We would like to use Theorem \ref{characterization}, for this we only need to check three conditions: equality of the kernels of $\theta$ and $\eta$, surjectivity of the representation $\theta\colon H_n(q)\ra\End_\Ccal(X^{\ten n})$ and $\Hom(X^{\ten m}, X^{\ten n})=\{0\}$ if $m\not\equiv n \pmod{d}$. Let us start with the easiest one: the last one. \\
For this note that $\Hom_\Ccal(X^{\ten m},X^{\ten n})$ is generated by $a$ and $\nu$. So if $\alpha\in\Hom_\Ccal(X^{\ten m},X^{\ten n})$, then $\alpha=0$, or $\alpha$ is a linear combination of words consisting of the letters $\iota^{\ten k}\ten \nu\ten\iota^{\ten l}$, $\iota^{\ten k}\ten\nu^*\ten\iota^{\ten l}$ and $\theta(x)$ for $k,l\in\Nat$ and $x\in H_\infty(q)$. It is sufficient to consider individual words. If $x\in H_p(q)\subset H_\infty(q)$, then $\theta(x)\in\End(X^{\ten p})$ and for $k,l\in\Nat$ we have $\iota^{\ten k}\ten \nu\ten\iota^{\ten l}\in\Hom(X^{k+l}, X^{k+d+l})$, $\iota^{\ten k}\ten \nu^*\ten\iota^{\ten l}\in\Hom(X^{k+d+l}, X^{k+l})$. Induction on the length of a word gives the result.\\
To be able to prove the other two remaining requirements we first compute $(\nu^*\ten\iota^{\ten k})(\iota^{\ten k}\ten\nu)$. 
In the upcoming computations we need the identity 
\begin{equation}\label{eq_embedding8}
\theta(g_i)\nu = \theta(g_i)\nu(\nu^*\nu) = \theta(g_iF_d)\nu = -\theta(F_d)\nu= -(\nu\nu^*)\nu= -\nu, \qquad \textrm{for } i=1,\ldots,d-1,
\end{equation}
which follows from \eqref{eq_embedding3}, \eqref{eq_embedding4} and Lemma \ref{properties_Bn}. This also implies that $\nu^*\theta(g_i)=-\nu^*$.

\begin{Lemma}\label{relations_nu_k}
For $k=1,2,\ldots,d-1$ the following equality holds
\begin{equation}\label{eq_relations_nu_k4}
(\nu^*\ten\iota^{\ten k})(\iota^{\ten k}\ten\nu)= \omega^k \,\frac{[d-k]_{\frac{1}{q}}![k]_{\frac{1}{q}}!}{[d]_{\frac{1}{q}}!} \, (-\mu)^{-k(d-k)}\theta(F_k).
\end{equation}
\end{Lemma}
\begin{pf}
We prove this by induction. The case $k=1$ is exactly assumption \eqref{eq_embedding5} of the theorem, so we will prove the induction step.
Consider the morphism $T:=(\nu^*\ten\iota^{\ten k })(\iota^{\ten k-1}\ten\nu\ten\iota)(\iota^{\ten k-1}\ten\nu^*\ten\iota)(\iota^{\ten k}\ten\nu)$. By the induction hypothesis and the assumption of this theorem, this morphism equals
\begin{align}
T&= \omega^{k-1} \,\frac{[d-k+1]_{\frac{1}{q}}![k-1]_{\frac{1}{q}}!}{[d]_{\frac{1}{q}}!}\, (-\mu)^{-(k-1)(d-k+1)} \omega (-\mu)^{-(d-1)}[d]_{\frac{1}{q}}^{-1}(\theta(F_{k-1})\ten\iota)\notag\\
&= \omega^{k} \,\frac{[d-k+1]_{\frac{1}{q}}![k-1]_{\frac{1}{q}}!}{[d-1]_{\frac{1}{q}}!}\, (-\mu)^{k^2 -kd-2k+2} (\theta(F_{k-1})\ten\iota).\label{eq_relations_nu_k2}
\end{align}
On the other hand as $\nu\nu^*=\theta(F_d)$, $\theta(g_i)\nu=-\nu$ and $\nu^*\theta(g_i)=-\nu^*$ we have 
\begin{align}
T&= (\nu^*\ten\iota^{\ten k})(\iota^{\ten k-1}\ten\theta(F_d)\ten\iota)(\iota^{\ten k}\ten\nu)\notag\\
&=[d]_{\frac{1}{q}}^{-1} (\nu^*\ten\iota^{\ten k}) \theta(1+ (-q)^{-1}g_k + \ldots + (-q)^{-(d-1)}g_{k+d-2} \cdots g_k)(\iota^{\ten k}\ten\theta(F_{d-1})\ten\iota) (\iota^{\ten k}\ten\nu)\notag\\
&=[d]_{\frac{1}{q}}^{-1} (\nu^*\ten\iota^{\ten k}) \theta(1+ q^{-1} + \ldots + q^{-(d-k)} + (-q)^{-(d+1-k)}g_d\cdots g_k \notag\\
&\qquad + \ldots + (-q)^{-(d-1)}g_{k+d-2} \cdots g_k)(\iota^{\ten k}\ten\nu)\notag \displaybreak[2]\\
&= [d]_{\frac{1}{q}}^{-1}[d+1-k]_{\frac{1}{q}} (\nu^*\ten\iota^{\ten k})(\iota^{\ten k}\ten\nu) + [d]_{\frac{1}{q}}^{-1} (-q)^{-(d+1-k)} \notag\\
&\qquad \cdot\theta(1+ (-q)^{-1}g_1+\ldots + (-q)^{-(k-2)}g_{k-2}\cdots g_1)(\nu^*\ten \iota^{\ten k})\theta(g_d\cdots g_k)(\iota^{\ten k}\ten\nu)\label{eq_relations_nu_k1}
\end{align}
Now note that by the assumptions and induction hypothesis
\begin{align*}
(\nu^*\ten &\iota^{\ten k})\theta(g_d\cdots g_k)(\iota^{\ten k}\ten\nu) = (\nu^*\ten \iota^{\ten k})\theta(g_{d+1}^{-1}\cdots g_{k+d-1}^{-1})\theta(g_{k+d-1}\cdots g_k)(\iota^{\ten k}\ten\nu)\\
&= \theta(g_1^{-1}\cdots g_{k-1}^{-1}) (\nu^*\ten \iota^{\ten k})\omega\mu^{d+1}(\iota^{\ten k-1}\ten\nu\ten\iota)\\
&= \omega\mu^{d+1} \omega^{k-1} \,\frac{[d-k+1]_{\frac{1}{q}}![k-1]_{\frac{1}{q}}!}{[d]_{\frac{1}{q}}!}\, (-\mu)^{-(k-1)(d-k+1)} \theta(g_1^{-1}\cdots g_{k-1}^{-1})(\theta(F_{k-1})\ten\iota).
\end{align*}
Since 
\[
(-q)^{-(d+1-k)}\mu^{d+1}(-\mu)^{-(k-1)(d-(k-1))}=(-1)^{k(d+1-k)}\mu^{-k(d-k)},
\]
identity \eqref{eq_relations_nu_k1} equals
\begin{align}
T=&[d]_{\frac{1}{q}}^{-1}[d+1-k]_{\frac{1}{q}} (\nu^*\ten\iota^{\ten k})(\iota^{\ten k}\ten\nu) + \omega^k [d]_{\frac{1}{q}}^{-1}(-1)^{k(d+1-k)}\mu^{-k(d-k)} \,\frac{[d-k+1]_{\frac{1}{q}}![k-1]_{\frac{1}{q}}!}{[d]_{\frac{1}{q}}!}\notag\\
&\qquad \cdot \theta(g_1^{-1}\cdots g_{k-1}^{-1} + (-q)^{-1}g_2^{-1}\cdots g_{k-1}^{-1}+\ldots + (-q)^{-(k-2)}g_{k-1}^{-1})(\theta(F_{k-1})\ten\iota).\label{eq_relations_nu_k3}
\end{align}
If we now combine both expressions of $T$, \eqref{eq_relations_nu_k2} and \eqref{eq_relations_nu_k3}, we get
\begin{align*}
&[d]_{\frac{1}{q}}^{-1}[d+1-k]_{\frac{1}{q}} (\nu^*\ten\iota^{\ten k})(\iota^{\ten k}\ten \nu) = (-1)^{k(d+1-k)+1} \mu^{-k(d-k)}\omega^k \,\frac{[d-k+1]_{\frac{1}{q}}![k-1]_{\frac{1}{q}}!}{[d]_{\frac{1}{q}}\,[d]_{\frac{1}{q}}!}\\
&\qquad\cdot \theta(g_1^{-1}\cdots g_{k-1}^{-1} + (-q)^{-1}g_2^{-1}\cdots g_{k-1}^{-1} + \ldots + (-q)^{-(k-1)}) (\theta(F_{k-1})\ten\iota)
\end{align*}
which by Lemma \ref{properties_Bn} equals
\[
(-1)^{k-1}(-1)^{k(d+1-k)+1} \mu^{-k(d-k)}\omega^k \,\frac{[d-k+1]_{\frac{1}{q}}![k-1]_{\frac{1}{q}}!}{[d]_{\frac{1}{q}}\,[d]_{\frac{1}{q}}!}\,[k]_{\frac{1}{q}}\theta(F_k).
\]
From this \eqref{eq_relations_nu_k4} follows immediately and the lemma is proved.
\end{pf}

\begin{Lemma}\label{equality_kernels}
The representation $\theta$ satisfies $\ker(\theta\colon H_n(q)\ra\End(X^{\ten n})) = \ker(\eta\colon H_n(q)\ra\End(\Hcal^{\ten n}))$. 
\end{Lemma}
\begin{pf}
From the above lemma it follows in particular that 
\[
(\nu^*\ten\iota^{\ten d-1})(\iota^{\ten d-1}\ten\nu) = \overline{\omega}(-\mu)^{-(d-1)}[d]_{\frac{1}{q}}^{-1}\,\theta(F_{d-1})
\]
and thus that the morphisms 
\[
R:= \overline{\omega}(-1)^{d-1}[d]_{\frac{1}{q}}^{1/2}\mu^{(d-1)/2}\nu, \qquad \overline{R}:=[d]_{\frac{1}{q}}^{1/2}\mu^{(d-1)/2}\nu
\]
satisfy the conjugate equations for $X$. Define a map 
\[
\varphi^{(n)}\colon \End(X^{\ten n})\ra \End(X^{\ten n-1}),\qquad \alpha\mapsto (\iota^{\ten n-1}\ten \nu^*)(\alpha\ten \iota^{\ten d-1})(\iota^{\ten n-1}\ten \nu)
\]
and the functional $\varphi_n:=\varphi^{(1)}\circ\cdots\circ\varphi^{(n-1)}\circ\varphi^{(n)}$. Now let $(R',\overline{R}')$ be a standard solution of the conjugate equations of $X$. The map 
\[
\End(X^{\ten n})\ra \End(X^{\ten n-1}),\qquad \alpha\mapsto (\iota^{\ten n-1}\ten \overline{R}'^*)(\alpha\ten \iota^{\ten d-1})(\iota^{\ten n-1}\ten \overline{R}') 
\]
is a partial trace induced by a standard solution, so it is tracial and faithful. There exists an invertible morphism $T\in\Hom(\overline{X}',\overline{X})$ such that $R=(T^{-1}\ten \iota)R'$ and $\overline{R}=(\iota\ten T^*)\overline{R}'$ \cite[Prop. 2.2.4]{NeshveyevTuset13}. From this it is immediate that $\varphi^{(n)}$ and thus $\varphi_n$ are also faithful. Using the involution, equation \eqref{eq_embedding5} can be rewritten as
\[
(\iota\ten\nu^*)(\nu\ten\iota)=\overline{\omega}(-\mu)^{-(d-1)}[d]_{\frac{1}{q}}^{-1}\iota.
\]
Combination with \eqref{eq_embedding7} and \eqref{eq_embedding8} gives that
\begin{align*}
\varphi^{(2)}\circ\theta(g_1) &= (\iota\ten\nu^*)\theta(g_1)(\iota\ten\nu)\label{eq_embedding9} \\
&= (-1)^{d-1}(\iota\ten\nu^*)\theta(g_d\cdots g_1)(\iota\ten\nu)\\
&= (-1)^{d-1}\omega\mu^{d+1}(\iota\ten\nu^*)(\nu\ten\iota)\\
&= (-1)^{d-1}\omega\mu^{d+1}\overline{\omega}(-\mu)^{-(d-1)}[d]_{\frac{1}{q}}^{-1}\iota\\
&= q\,[d]_{\frac{1}{q}}^{-1}\iota = \frac{q^d}{[d]_q}\,\iota.
\end{align*}
Thus in particular $\varphi^{(2)}\circ\theta(g_1)$ is a scalar in $\End(X)$ and thus $\varphi^{(2)}\circ\theta(e_1)$ is a scalar as well. Therefore if $x,y\in H_{n-1}(q)$  
\begin{align*}
\varphi^{(n)}(\theta(xe_{n-1}y))&= (\iota^{\ten n-1}\ten \nu^*)\theta(x)(\iota^{\ten n-2}\ten \theta(e_1))\theta(y)(\iota^{\ten n-1}\ten \nu)\\
&= \theta(x)(\iota^{\ten n-1}\ten \nu^*)(\iota^{\ten n-2}\ten \theta(e_1))(\iota^{\ten n-1}\ten \nu)\theta(y)\\
&= \varphi_2(\theta(e_1))\cdot\theta(xy).
\end{align*}
So $\varphi_n\circ\theta$ defines a faithful functional with the Markov property on $H_n(q)$. According to Lemma \ref{equivalent_Markov_traces} this functional must be tracial and hence we obtain a Markov trace $\tr_\Ccal := \varphi_n\circ\theta\colon H_n(q)\ra\Com$. Markov traces are characterized by their value on the generator $g_1$. Recall from Theorem \ref{equivalence_of_representations} that $\tr_{\Rep(SU_\mu(d))}(g_1)= \frac{q^d}{[d]_q}$. It follows that $\tr_\Ccal = \tr_{\Rep(SU_\mu(d))}$ and thus $\ker(\theta)=\ker(\tr_\Ccal)=\ker(\tr_{\Rep(SU_\mu(d))})=\ker(\eta)$.
\end{pf}

Now surjectivity of the map $\theta\colon H_n(q)\ra\End(X^{\ten n})$. For this we need the following lemma. Recall the notation $\nu_{k,l}:=\iota^{\ten k}\ten \nu\ten\iota^{\ten l}$ and $\nu^*_{k,l}:=\iota^{\ten k}\ten \nu^*\ten\iota^{\ten l}$.

\begin{Lemma}\label{reduction}
Let $x\in H_p(q)$ and $k,l,m,n,k',l',m',n'\in\Nat$ be natural numbers satisfying the equality $k+l+d=m+n+d = k'+l'=m'+n'= p$, then there exist $x_1,x_2\in H_\infty(q)$ such that 
\[
\nu_{k,l}^*\theta(x)\nu_{m,n} = \theta(x_1),\qquad \nu_{k',l'}\theta(x)\nu_{m',n'}^*=\theta(x_2).
\]
\end{Lemma}
\begin{pf}
First we prove the second assertion. We write $k,l,m,n$ instead of $k',l',m',n'$. Note that by \eqref{eq_embedding7} there exist $y_1\in H_{k+d}(q)$ and $y_2\in H_{m+d}(q)$ such that 
\begin{equation}\label{eq_reduction1}
\nu_{k,l}=\theta(y_1)\nu_{0,k+l},\qquad \nu_{m,n}^*=\nu_{0,m+n}^*\theta(y_2).
\end{equation}
Then by \eqref{eq_embedding4}
\begin{align*}
\nu_{k,l}\theta(x)\nu_{m,n}^* &= \theta(y_1)\nu_{0,k+l}\theta(x)\nu_{0,m+n}^*\theta(y_2) = \theta(y_1) (\iota^{\ten d}\ten\theta(x))(\nu\nu^*\ten\iota_{k+l})\theta(y_2)\\
& = \theta(y_1\Sigma^d(x)F_dy_2),
\end{align*}
where we still use $\Sigma$ to denote the shift map. \\ 
Now the first case. Similar to \eqref{eq_reduction1} there exist $y_1$ and $y_2$ such that 
\[
\nu_{k,l}^*\theta(x)\nu_{m,n} = \nu_{0,k+l}\theta(y_1)\theta(x)\theta(y_2)\nu^*_{0,m+n}.
\]
So we can assume to deal with the case $\nu_{0,k}^*\theta(x)\nu_{0,k}$ and $x\in H_{d+k}(q)$. Now observe that $\nu_{0,k}^*\theta(x)\nu_{0,k} = \nu_{0,k}^*\theta(F_dxF_d)\nu_{0,k}$. By surjectivity of the representation $\eta\colon H_k(q)\ra\End(\Hcal^{\ten k})$ there exists an $y\in H_k(q)$ such that $(S^*\ten\iota^{\ten k})\eta(F_dxF_d)(S\ten\iota^{\ten k})=\eta(y)$, here $S$ and $\eta$ are as in Notation \ref{Def_rep_SUmud}. This implies that 
\[
\eta(F_d\Sigma^d(y)) = SS^*\ten\eta(y) = (S\ten\iota^{\ten k})\eta(y)(S^*\ten\iota^{\ten k}) = (SS^*\ten\iota^{\ten k})\eta(F_d x F_d) (SS^*\ten\iota^{\ten k}) = \eta(F_d x F_d).
\]
Because by Lemma \ref{equality_kernels} the representations $\eta$ and $\theta$ have the same kernel, it follows that $\theta(F_dxF_d)=\theta(F_d\Sigma^d(y))$. Combining all this gives
\[
\nu_{0,k}^*\theta(x)\nu_{0,k}= \nu_{0,k}^*\theta(F_dxF_d)\nu_{0,k} = \nu_{0,k}^*\theta(F_d\Sigma^d(y))\nu_{0,k} = \theta(y)\nu_{0,k}^*\nu_{k,0}\nu_{0,k}^*\nu_{k,0} = \theta(y)
\]
and concludes the lemma.
\end{pf}

To prove that the representation $\theta\colon H_n(q)\ra\End_\Ccal(X^{\ten n})$ is surjective let $\alpha\in\End_\Ccal(X^{\ten n})$. Then $\alpha$ is a linear combination of words consisting of the letters $\theta(x)$ for $x\in H_\infty(q)$ and $\nu_{k,l}$, $\nu_{m,n}$. Let $\beta=\beta_1\cdots\beta_r$ be such a word and $\beta_i$ the letters. Then $\beta\in\End_\Ccal(X^{\ten n})$ and thus 
\[
\#\{i\,:\, \beta_i=\nu_{k,l}\textrm{ some } k,l\} = \#\{i\,:\, \beta_i=\nu^*_{k,l}\textrm{ some } k,l\}.
\]
We now apply induction on $r$. If $r=1$, then the above sets must be empty and thus $\beta=\theta(x)$ for some $x\in H_n(q)$. Suppose $r>1$ and not all $\beta_i$ are of the form $\theta(x)$ for $x\in H_\infty(q)$, then there must exist $1\leq i<j\leq r$ such that either $\beta_i=\nu_{k,l}$, $\beta_j=\nu_{m,n}^*$ for some $k,l,m,n$ and $\beta_s=\theta(x_s)$ for all $i<s<j$, $x_s\in H_\infty(q)$ or $\beta_i=\nu^*_{k,l}$, $\beta_j=\nu_{m,n}$ for some $k,l,m,n$ and $\beta_s=\theta(x_s)$ for all $i<s<j$, $x_s\in H_\infty(q)$. In both cases we can apply Lemma \ref{reduction} to reduce $\beta_i\beta_{i+1}\cdots\beta_j$ to $\theta(x)$ for some $x\in H_\infty(q)$. In this way we obtain a word of length $<r$ and by induction $\beta\in \theta(H_n(q))$. Hence $\theta\colon H_n(q)\ra\End_\Ccal(X^{\ten n})$ is surjective and the conclusion follows from Theorem \ref{characterization}.
\end{pf}

\section*{Acknowledgements}
The author would like to thank Sergey Neshveyev for his help and for carefully reading the manuscript. 

\bibliographystyle{hplain}
\bibliography{/mn/anatu/gjester-u2/bpjordan/Documents/Papers/References/references}
\end{document}